\documentclass[twocolumn]{autart_modified}    

\usepackage[dvipdfmx]{graphics} 
\usepackage[cmex10]{amsmath} 
\usepackage{amssymb}  

\usepackage{cases}

\usepackage[subrefformat=parens]{subcaption} 

\usepackage{graphicx} 

\usepackage[dvipdfmx,usenames]{color}

\usepackage{cite}

\usepackage{algorithm}
\usepackage{algorithmic}

\usepackage[overlay]{textpos}

\usepackage{pgfplots}
\pgfplotsset{compat=newest}
\usepackage{tikz}
\usetikzlibrary{positioning}
\usetikzlibrary{arrows}

\newcommand{\argmax}{\operatornamewithlimits{\mathrm{arg~max}}}
\newcommand{\argmin}{\operatornamewithlimits{\mathrm{arg~min}}}

\begin{document}

\begin{frontmatter}

\title{Pseudo-perturbation-based Broadcast Control \\ of Multi-agent Systems\thanksref{footnoteinfo}} 

\thanks[footnoteinfo]{%
The material herein was presented in part at the 2018 American Control Conference \cite{ItoPBC2017}. 
Corresponding author is Yuji Ito, Tel. +81-561-71-7512.
Md Abdus Samad Kamal is currently with the School of Engineering, Monash University Malaysia, Bandar Sunway, 47500 Selangor Darul Ehsan, Malaysia.
Shun-ichi Azuma is currently with the Graduate School of Engineering, Nagoya University, Furo-cho, Chikusa-ku, Nagoya-shi, Aichi 464-8603, Japan.
} 

\author[tytlabs]{Yuji Ito}\ead{ito-yuji@mosk.tytlabs.co.jp},    
\author[tytlabs]{Md Abdus Samad Kamal}\ead{maskamal@ieee.org},               
\author[tytlabs]{Takayoshi Yoshimura}\ead{yoshimura@mosk.tytlabs.co.jp},  
\author[Kyoto]{Shun-ichi Azuma}\ead{shunichi.azuma@mae.nagoya-u.ac.jp}  

\address[tytlabs]{TOYOTA CENTRAL R\&D LABS., INC., Nagakute-shi, Aichi 480-1192, Japan}  
\address[Kyoto]{Graduate School of Informatics, Kyoto University, Kyoto-shi, Kyoto 606-8501, Japan}        

\begin{keyword}         
Broadcast control; Multi-agent systems; Simultaneous perturbation stochastic approximation; Stochastic control.        
\end{keyword}

\begin{abstract}                        
The present paper proposes a novel broadcast control (BC) law for multi-agent coordination. 
A BC framework has been developed to achieve global coordination tasks with low communication volume.
The BC law uses broadcast communication, which transmits an identical signal to all agents indiscriminately without any agent-to-agent communication.
Unfortunately, all of the agents are required to take numerous random actions because the BC law is based on stochastic optimization.
Such random actions degrade the control performance for coordination tasks and may invoke dangerous situations.
In order to overcome these drawbacks, the present paper proposes the pseudo-perturbation-based broadcast control (PBC) law, which introduces multiple virtual random actions instead of the single physical action of the BC law. 
The following advantages of the PBC law are theoretically proven.
The PBC law achieves coordination tasks asymptotically with probability 1.
Compared with the BC law, unavailing actions are reduced and agents' states converge at least twice as fast.
Increasing the number of multiple actions further improves the control performance because averaging multiple actions reduces unavailing randomness. 
Numerical simulations demonstrate that the PBC improves the control performance as compared with the BC law.

\end{abstract}

\end{frontmatter}

\section{Introduction} \label{sec_intro}

Coordinations of multi-agent systems have recently attracted attention in various engineering fields in line with recent technological advancements in computing and communication resources. 
Swarm robots, groups of vehicles, embedded robotic systems, and electric power networks are examples of multi-agent systems that require coordination in order to achieve desired global behaviors. 
Some common behaviors desired in multi-agent systems include coverage control by agents deployed over a given space, rendezvous with formation of the agents in a synchronized manner, and assignment at target locations (for examples, see \cite{martinez2007motion, Ren11, Ji06, kamal2007coordination}).

Decentralized/distributed control approaches perform multi-agent tasks using agent/supervisor-to-agent communication, e.g., \cite{Mahmoud11,ji2007rendezvous}.
In distributed frameworks, each agent takes action locally and individually, communicating with its neighbor agents. 
Alternatively, centralized control systems use a supervisor to aggregate information from all agents.
Such aggregation is efficient for achieving \textit{global} coordination tasks.
The supervisor evaluates a global objective based on information of agents and transmits commands to all of the agents individually. 
However, such centralized systems require a massive communication volume between a supervisor and agents, i.e., \textit{one-to-all} and \textit{all-to-one} communications.

For achieving global coordination tasks with a low communication volume, the present paper focuses on another structure of multi-agent systems called a broadcast system.
The broadcast system can work with a low communication volume because a supervisor broadcasts an identical signal to all of the agents indiscriminately and no agent-to-agent communication is used.
Table \ref{tab:communications} indicates that the both-way communication volume ($O(N+1)$) of the broadcast system is lower than that ($O(2N)$) of centralized systems employing unicast protocols in one-to-all communication.
The broadcast system has the potential to realize global behaviors because the supervisor aggregates the information of all of the agents.
Broadcast systems are constructed even if supervisors and all agents use broadcast communication rather than unicast communication, e.g., connected vehicles in traffic coordination \cite{KishimotoARIB2014}.
Combining broadcast communication with agent-to-agent communication ensures controllability and/or stability of distributed systems \cite{yoon2014controllability, yoon2014broadcasting}.

%

\begin{table}
	\renewcommand{\arraystretch}{0.2}	
	\caption{Communication volumes of unicast and broadcast protocols. The number of agents is denoted by $N$.}
	\label{tab:communications}
	{\scriptsize 
		\begin{center}
			\begin{tabular}{| c|| c |c | c | }
				\hline
				Protocol  
				&  \shortstack{one-to-all \\ communication}  &  \shortstack{all-to-one \\ communication} \\ \hline	 
				Unicast 		  &  $O(N)$ &	$O(N)$	\\ \hline	 
				Broadcast	  &  $O(1)$ &	$O(N)$	\\ \hline	 
			\end{tabular}	
		\end{center}
	}
	\vspace*{-0.0in}
\end{table}

The concept of broadcast control (BC) \cite{azuma2013broadcast} has recently been proposed for multi-agent coordination, which is applied to broadcast systems.
The concept is related to stochastic/biological systems \cite{ueda2007broadcast, halasz2008stochastic}.
The BC law can achieve global objectives, such as coverage, assembly, and assignment tasks \cite{Mohamad16assembly,Mohamad16assignment}.
A key scheme of the BC is inspired by simultaneous perturbation stochastic approximation (SPSA) \cite{Spall94,Spall98,Ahmad15}.
After each agent randomly moves temporarily, a supervisor broadcasts a signal that corresponds to the achievement degree of an objective, i.e., the value of the corresponding objective function.
The agents determine their next actions from the received signal and their random movements.
These actions implicitly construct an approximate gradient of the objective function.
The BC law employs a stochastic gradient method to minimize the objective function in the manner of SPSA.
However, such a process restricts its applicability. 
The random actions of agents involve significant traveling cost and time in order to achieve an objective. 
In many systems, e.g., vehicular traffic systems, taking random actions should be avoided in consideration of safety and operational efficiency.

In order to overcome the limitations of the BC law, the present paper proposes a control law called the pseudo-perturbation-based broadcast control (PBC) law.
In order to avoid unavailing random actions and improve the control performance for coordination tasks, a key concept is employing \textit{multiple virtual} random actions instead of the single physical random action of the BC.
The PBC for broadcast systems is operated in a manner similar to the BC and retains the basic features of the BC.
The present paper theoretically proves the following advantages and characteristics of the PBC law.
Coordination tasks for broadcast systems are asymptotically achieved with probability $1$.
Unavailing random actions of agents are reduced and coordination tasks are achieved at least twice as fast, as compared to the BC law.
The execution of a task is accelerated by increasing the number of multiple actions because the multiple actions improve the approximation accuracy of the gradient of an objective function for the task.
The PBC law is demonstrated through numerical simulations of two types of multi-agent coordination tasks, which are a coverage task and a rendezvous task with formation selection.

This paper is an extension of our previous work (under review) \cite{ItoPBC2017}.
We modified the paper significantly and added new points.
Especially, this paper proposes the PBC law based on \textit{multiple} virtual actions of agents, whereas the control law in the conference paper uses a \textit{single} virtual action. 
To discuss \textit{multiple} actions, the overview of the PBC in Section \ref{sec_PBC_overview} and Theorem \ref{thm:convergence_of_PBC} were modified. 
All of Section \ref{sec_PBC_analysis2} is a new part. 
Numerical simulations in Sections \ref{sec_sim_rendezvous} and \ref{sec_sim_coverage} to demonstrate the PBC law are new.
All appendices to prove theorems and a proposition are new, where Appendices \ref{pf:PBC_cost_improvement} and \ref{pf:PBC_dist_reduction} were a modified version of \cite{ItoPBC2017}.

The remainder of this paper is organized as follows.
Section \ref{sec_problem} states the primary problem and associated limitations of the BC law. 
The PBC law with theoretical analysis is presented in Section \ref{sec_PBC}, which describes the main results of the paper.
In order to demonstrate the effectiveness of the PBC law, two types of multi-agent coordination tasks are simulated in Section \ref{sec_simulation}.
Conclusions and future research are described in Section \ref{sec_conclusion}.

\textbf{Notation:} 
For a vector $x:=[x_{1},...,x_{n}]^{\mathrm{T}} \in \mathbb{R}^{n}$ with nonzero elements, $x^{(-1)}$ represents the element-wise inverse of $x$, i.e.,  $x^{(-1)}:=[1/x_{1},...,1/x_{n}]^{\mathrm{T}} \in \mathbb{R}^{n}$. 
The partial derivative of $f(x) \in \mathbb{R}$ with respect to $x \in \mathbb{R}^{n}$ is denoted by $\partial_{x} f(x) \in \mathbb{R}^{n}$.
The operator $\mathrm{E}[f(\sigma)]$ denotes the expectation of $f(\sigma) \in \mathbb{R}^{n}$ with respect to all random variables $\sigma$.
The operators $\mathrm{E}[f(\sigma,x(\sigma))|_{x}]$ and $\mathrm{Cov}[f(\sigma,x(\sigma))|_{x}]$ denotes the conditional expectation and the conditional covariance of $f(\sigma,x(\sigma)) \in \mathbb{R}^{n}$ with respect to all random variables $\sigma$ for a fixed $x$, respectively.


\newcommand{\numK}{K}

\newcommand{\defGrad}{g}
\newcommand{\ERRdefGrad}{e}

\newcommand{\BCnoiseRand}{e_{\mathrm{M}}}
\newcommand{\BCnoiseOffset}{e_{\mathrm{B}}}

\newcommand{\totaldist}{D}

\newcommand{\AssignID}{S}

\newcommand{\SmallParamCtwo}{\varepsilon}

\newcommand{\JofGrad}{J_{\defGrad}}
\newcommand{\valFunc}{V}

\section{Problem setting via review of the BC} \label{sec_problem}

This section states the main problem and associated drawbacks of the BC law.
Section \ref{sec_target} introduces the target systems and the main problem.
After Section \ref{sec_review_BC} reviews the BC law, which is a solution to the main problem, its limitations to be overcome are described in Section \ref{sec_main_problem}.

%
%

\subsection{Target systems and the main problem} \label{sec_target}

Let us consider a multi-agent system $\Sigma$ on an $n$-dimensional state space, which is controlled as a broadcast system.
The system consists of a global controller (supervisor) $G$ and $N$ agents $A_{i}$ equipped with local controllers $L_{i}$ ($i = 1,...,N$).
Note that, although agents are indexed from $1$ to $N$ for notational convenience, it is not necessary to discriminate among them.

The dynamics of the agent $A_{i}$ at a discrete time step $t \in \{0,1,...\}$ is omni-directional and is given by
\begin{equation}\label{eq:Agent}
	\begin{aligned}
		A_{i}: x_{i}(t+1)=x_{i}(t)+u_{i}(t)
		,
	\end{aligned}
\end{equation}
where $x_{i}(t) \in \mathbb{R}^n$ and $u_{i}(t) \in \mathbb{R}^n$ are the state and the control input of $A_{i}$, respectively.
The collective state of all of the agents is denoted by $x:=[x_{1}^{\mathrm{T}},..., x_{N}^{\mathrm{T}}]^{\mathrm{T}}\in \mathbb{R}^{nN}$.
Its initial state $x(0)$ is given.
The local controller $L_{i}$, with which each agent $A_{i}$ is equipped, is described as follows:
\begin{equation}\label{eq:LocalC}
\begin{aligned}
L_{i}: 
\begin{cases}
\phi_{i}(t+1) &= {f}_{\phi}(\phi_{i}(t),\nu(t),t), \\
\qquad u_{i}(t) &=  {f}_{u}(\phi_{i}(t),\nu(t),t),
\end{cases}
\end{aligned}
\end{equation}
where $\phi_{i} \in \mathbb{R}^{n_{\phi}}$ is the state of the local controller $L_{i}$.
The collective state of all of the local controllers is denoted by $\phi:=[\phi_{1}^{\mathrm{T}},...,\phi_{N}^{\mathrm{T}}]^{\mathrm{T}} \in \mathbb{R}^{n_{\phi}N}$.
The initial state is set to $\phi(0)=0 \in \mathbb{R}^{n_{\phi}N}$ for simplicity of discussion.
The symbol $\nu(t) \in \mathbb{R}^{n_{\nu}}$ is a broadcast signal, which is sent to all of the agents by the global controller $G$. 
The functions ${f}_{\phi}:\mathbb{R}^{n_{\phi}} \times \mathbb{R}^{n_{\nu}} \times \mathbb{Z} \rightarrow \mathbb{R}^{n_{\phi}}$ and ${f}_{u}:\mathbb{R}^{n_{\phi}} \times \mathbb{R}^{n_{\nu}} \times \mathbb{Z}  \rightarrow \mathbb{R}^{n}$ determine the transition of $\phi_{i}$ and the input $u_{i}(t)$ of the agent $A_{i}$, respectively.
The global controller $G$ is expressed as
\begin{equation}\label{eq:globalC}
\begin{aligned}
G: \nu(t)=f_{\nu}(x(t), \phi(t),t) \in  \mathbb{R}^{n_{\nu}}
,
\end{aligned}
\end{equation}
where ${f}_{\nu}:\mathbb{R}^{nN} \times \mathbb{R}^{n_{\phi}N}  \times \mathbb{Z} \rightarrow \mathbb{R}^{n_{\nu}}$ determines the broadcast signal $\nu(t)$.

Note the following three requirements with respect to the communication and controllers in terms of cost and equipment for the multi-agent system $\Sigma$.
First, there is no agent-to-agent communication.
Second, the local controllers ${f}_{\phi}$ and ${f}_{u}$ must be identical for all of the agents. 
Third, the broadcast signal $\nu(t)$ must be identical for all of the agents (each optimal signal cannot be sent to each agent).
Standard centralized, decentralized, and distributed control approaches can no longer be applied due to these requirements, whereas the requirements are crucial for designing controllers for large-scale multi-agent systems.
The main problem of realizing coordination tasks under these requirements is stated as follows.

\textbf{Main problem:}
For the multi-agent system $\Sigma$, find global and local controllers $(G,L_{1},L_{2},\ldots,L_{N})$ (i.e., $f_{\nu}$, ${f}_{\phi}$, and ${f}_{u}$) that satisfy the three above-mentioned requirements and
\begin{equation}\label{eq:LimJ}
\begin{aligned}
\lim_{t\rightarrow \infty} J(x(t))=\min_{x\in \mathbb{R}^{nN}}J(x)
,
\end{aligned}
\end{equation}
where $J: \mathbb{R}^{nN}\rightarrow \mathbb{R}$ is an objective function that is regarded as the performance index for a coordination task of the system $\Sigma$.
Convergence to a local minimum of $J(x)$ in (\ref{eq:LimJ}) is allowed in this problem.

\subsection{Review of the BC law} \label{sec_review_BC}

This subsection briefly introduces the BC law \cite{azuma2013broadcast}.
The global controller $G$ of the BC law is defined as
\begin{equation}\label{eq:G_BC}
\begin{aligned}
\nu(t)=f_{\nu}(x(t), \phi(t),t):=J(x(t))
\in 
\mathbb{R}
.
\end{aligned}
\end{equation}
Each local controller $L_{i}$ is defined by
\begin{align}
\phi_{i}(t) &:= 
\begin{bmatrix}
\phi_{i,1}(t) \\
\phi_{i,2}(t) 
\end{bmatrix}
\in
\{-1,1\}^{n} \times \mathbb{R}
,\label{eq:L_state_BC}
\\
\phi_{i}(t+1)
&\;=
{f}_{\phi}(\phi_{i}(t),\nu(t),t) := 
\begin{bmatrix}
\sigma_{i}(t) \\
\nu(t)
\end{bmatrix}
,\label{eq:L1_BC}
\\
u_{i}(t)&\;={f}_{u}(\phi_{i}(t),\nu(t),t) 
\nonumber\\&
:=
\begin{cases}
c(t)\sigma_{i}(t)  \qquad ( t \in \{0,2,4,...\} ) \\
-c(t)\phi_{i,1}(t)-a(t)\frac{\nu(t)-\phi_{i,2}(t)}{c(t)} \phi_{i,1}^{(-1)}(t) \\
\qquad\qquad\quad\;\;  ( t \in \{1,3,5,...\} )
,
\end{cases}
\label{eq:L2_BC}
\end{align}
where $\phi_{i,1}(t) \in \{-1,1\}^n$ and $\phi_{i,2}(t) \in \mathbb{R}$.
The symbols $a(t) \in \mathbb{R}$ and $c(t)  \in \mathbb{R}$ are controller gains.
The symbol $\sigma_{i}(t):=[\sigma_{i,1}(t),...,\sigma_{i,n}(t)]^{\mathrm{T}} \in \{-1,1\}^{n}$ is a random variable, where $\sigma_{i,j}(t)$ ($i \in \{1,...,N\}$, $j \in \{1,...,n\}$, $t \in \{0,1,...\}$) independently obeys the Bernoulli distribution with outcome $\pm 1$ equal probabilities.
The collection of the random variables is denoted by $\sigma(t):=[\sigma_{1}(t)^{\mathrm{T}},...,\sigma_{N}(t)^{\mathrm{T}}]^{\mathrm{T}} \in \{-1,1\}^{nN}$.
The BC law with these definitions is a solution to the main problem.

\begin{thm}[BC law \cite{azuma2013broadcast}] \label{thm:BC}
For the multi-agent system $\Sigma$, an objective function $J(x)$, and the BC law in (\ref{eq:G_BC}), (\ref{eq:L_state_BC}), (\ref{eq:L1_BC}), and (\ref{eq:L2_BC}), $x(t)$ converges to a (possibly sample-path-dependent) solution to $\partial_{x}J(x) =0$ with probability $1$ if the following conditions hold.
\begin{enumerate}
	\item [(c1)] 
	An objective function  $J: \mathbb{R}^{nN} \rightarrow \mathbb{R}$ is defined as 
	\begin{align}
	J(x(t))&:=\rho(\|x\|)J_{\mathrm{obj}}(x)  + ( 1 - \rho(\|x\|) ) x^{\mathrm{T}}x
	,\label{eq:barrier_func}
	\\
	\rho(\|x\|)&:=
	\begin{cases}
	1 & ( \|x\| \leq l_{1}),
	\\
	0 & ( \|x\| \geq l_{2}),
	\end{cases}
	\label{eq:barrier_switch}
	\end{align}
	where $J(x)$ and $J_{\mathrm{obj}}(x)$ are nonnegative $C^{2}$ continuous on $\mathbb{R}^{nN}$, and there exists a solution to $\partial_{x}J(x) =0$ for some $l_{1}$ and $l_{2}$.	
	\item [(c2)]
	The compact connected internally chain transitive invariant sets%
	\footnote{%
		Consider the system $\dot{z}(\tau)=f(z(\tau))$ and one of its invariant sets, $\mathbb{S}_{z}$. 
		The invariant set $\mathbb{S}_{z}$ is said to be internally chain transitive if, for each $(z_{0}, z_{f}) \in \mathbb{S}_{z}\times\mathbb{S}_{z}$, $\epsilon > 0 \in \mathbb{R}$, and $T > 0 \in \mathbb{R}$, there exist $m \in \mathrm{N}$ and $(z_{1},z_{2},...,z_{m}) \in \mathbb{S}_{z}^{m}$ such that $\|z(\tau,z_{i})-z_{i+1}\|< \epsilon$ $(i = 0, 1,..., m)$ for some $\tau \in [T, \infty)$, where $z_{m+1}:=z_{f}$ and $z(\tau, z_{i})$ represents the state $z(\tau)$ of the system for the initial state $z(0) := z_{i}$ \cite{Borkar08}.
	}%
	of a gradient system $\dot{z}(\tau)=-\partial_{z}J(z(\tau))$ are included in the solution set to the equation 
	$\partial_{z}J(z)=0$ (i.e., to $\partial_{x}J(x)=0$),
	and there exists an asymptotically stable equilibrium for the gradient system, where $z(\tau) \in \mathbb{R}^{nN}$ and the stability is in the Lyapunov sense.	
	\item [(c3)]
	The controller gains satisfy $a(t)=a(t+1) > 0$ and $c(t)=c(t+1) > 0$ for every $t \in \{0, 2, 4,...\}$ (i.e., $a(0) = a(1) > 0$, $a(2) = a(3) > 0$,... and the same is true of $c(t)$), $\lim_{t \rightarrow \infty}a(t)=0$, $\sum_{t=0}^{\infty}a(t)=\infty$, $\lim_{t \rightarrow \infty}c(t)=0$, $\sum_{t=0}^{\infty}(a(t)/c(t))^{2}<\infty$, and $\sum_{t=0}^{\infty}a(t)c(t)^{2}<\infty$.	
\end{enumerate} 
\end{thm}
\begin{rem}\label{rem:BC_solution}
Convergence to a solution to $\partial_{x}J(x) =0$ indicates that (\ref{eq:LimJ}) (approximately) holds.
Thus, the BC law is a solution to the main problem.
\end{rem}
\begin{rem}\label{rem:BC_explanation}
	The BC law is based on a two-stage transition of the state $x(t)$.
	At the time $t \in \{0,2,4,...\}$, the agent $A_i$ receives the broadcast signal $\nu(t)$ and takes the random action $u_{i}(t)=c(t)\sigma_{i}(t)$. 
	Such a random action is essential in order to approximate the gradient of $J(x)$ based on SPSA \cite{Spall92}.
	The second term in (\ref{eq:L2_BC}) corresponds to the approximate gradient (in an average sense) by applying Taylor's theorem \cite{azuma2013broadcast}
	\begin{equation}\label{eq:E_grad_BC}
		\begin{aligned}
		&
		\forall t \in \{ 0,2,4,... \} , 
		\\&
		\mathrm{E}\Big[
		\frac{\nu(t+1)-\phi_{i,2}(t+1)}{c(t+1)} \phi_{i,1}^{(-1)}(t+1)
		\Big|_{x(t)}
		\Big]	
		\\&
		=
		\mathrm{E}\Big[
		\frac{J(x(t)+c(t)\sigma(t))-J(x(t))}{c(t)} \sigma_{i}^{(-1)}(t)
		\Big|_{x(t)}
		\Big]
		\\&
		=
		\partial_{x_{i}}J(x(t)) + O(c(t)) 
		,\; \mathrm{as} \;
		c(t) \rightarrow 0
		,	
		\end{aligned}
	\end{equation}
	where
	$\nu(t+1)=J(x(t+1))$ from (\ref{eq:G_BC}), 
	$x_{i}(t+1)=x_{i}(t)+c(t)\sigma_{i}(t)$ from (\ref{eq:Agent}) and (\ref{eq:L2_BC}),  
	$[\phi_{i,1}(t+1)^{\mathrm{T}}, \phi_{i,2}(t+1)] = [\sigma_{i}(t)^{\mathrm{T}}, J(x(t))]$ from (\ref{eq:L_state_BC}) and (\ref{eq:L1_BC}),
	and
	$c(t+1)=c(t)$ from the condition (c3)
	 hold.
	The gradient becomes the moving direction of $A_i$ for minimizing $J(x)$.  
	At the time $t+1 \in \{1,3,5,...\}$, $A_i$ moves in such a direction while canceling the random action $c(t+1)\phi_{i,1}(t+1)=c(t)\sigma_{i}(t)$. 		
\end{rem}
\begin{rem}
In the condition (c1), $J_{\mathrm{obj}}: \mathbb{R}^{nN} \rightarrow \mathbb{R}$ describes the achievement degree of a coordination task to be realized, such as coverage, rendezvous, and assignment tasks.
The function $\rho(\|x\|)$ switches $J(x)$ such that $J(x)=J_{\mathrm{obj}}(x)$ for $\|x\| \leq l_{1}$ and $J(x)=x^{\mathrm{T}}x$ for $\|x\| \geq l_{2}$.
The term $x^{\mathrm{T}}x$ is introduced for bounding the size of $x$.
The constants $l_{1}$ and $l_{2}$ determine the workspace of the coordination.
There exists $\rho(\|x\|)$ in (\ref{eq:barrier_switch}) such that $J(x)$ is $C^{2}$ continuous, e.g., a spline type function \cite{azuma2013broadcast}.		
\end{rem}
\begin{rem}\label{rem:BC_com_cost}
To evaluate $J(x(t))$ (if it is not directly observed), the global controller $G$ requires the information of all agents' states $x(t) \in \mathbb{R}^{nN}$, via some means, such as cameras, sensors, and direct \textit{all-to-one} communications.
\end{rem}

\subsection{Drawbacks of the BC law} \label{sec_main_problem}

For large-scale multi-agent systems, the BC law in Theorem \ref{thm:BC} is advantageous because it satisfies the three requirements described in Section \ref{sec_review_BC}: communications between agents are not used, local controllers are identical for all agents, and the broadcast signal is identical.
Such characteristics are efficient for constructing controllers that are not adversely affected by the system size.
However, the BC law contains the following limitations.
\begin{itemize}
	\item 
	The control law using the random action $c(t)\sigma_{i}(i)$ in (\ref{eq:L2_BC}) causes an agent to move an unavailing extra distance to approximate the gradient of the objective function $J(x)$. 
	If a broadcast system involves a disturbance (e.g., noise in observing the states of the agents), a large random action may be needed for the system to be robust against the disturbance.
	Furthermore, it is impractical to enforce random actions of agents in many applications because such a randomness invokes dangerous situations.
	\item 
	The minimization process in (\ref{eq:LimJ}) using the BC law can be slow because the state transition is a two-stage transition with randomness, as shown in (\ref{eq:L2_BC}). 
\end{itemize}
In order to overcome the above limitations, the next section presents a novel control law, which is another solution to the main problem, such that the moving distances of agents are reduced and the minimization of $J(x)$ is accelerated as compared with the BC law.

\section{Proposed method: pseudo-perturbation-based broadcast control law} \label{sec_PBC}

In this section, we propose the PBC law in order to overcome the drawbacks of the BC law.
Section \ref{sec_PBC_overview} presents an overview of the PBC law based on a simple but novel concept.
Sections \ref{sec_PBC_analysis1} and \ref{sec_PBC_analysis2} derive the theoretical aspects that show the effectiveness of the PBC law as compared with the BC law.

\subsection{Overview of the PBC law based on a simple but novel concept} \label{sec_PBC_overview}

The PBC law is derived by extending the BC law as follows.
Recall that the agents under the BC law take the random actions $u_{i}(t)=c(t)\sigma_{i}(t)$ in (\ref{eq:L2_BC}) in order to approximate the gradient of an objective function $J(x)$, as explained in Remark \ref{rem:BC_explanation}.
In order to avoid the unavailing random actions and improve the approximation accuracy of the gradient, a key concept is the introduction of the following \textit{multiple virtual} random actions rather than the single physical random actions:
\begin{align}
\hat{x}_{i}^{(k)}(t+1)&:=x_{i}(t)+\hat{u}_{i}^{(k)}(t) 
,\;&(k=1,2,...,\numK)
,\label{eq:x_predictive}
\\
\hat{u}_{i}^{(k)}(t)&:=c(t)\sigma_{i}^{(k)}(t)
,\;&(k=1,2,...,\numK)
,\label{eq:u_predictive}
\end{align}
where $\hat{u}_{i}^{(k)}(t)$ and $\hat{x}_{i}^{(k)}(t+1)$ are a virtual input and a virtual predictive state of each agent $A_{i}$, respectively.
The global controller $G$ calculates $\hat{u}_{i}^{(k)}(t)$ and $\hat{x}_{i}^{(k)}(t+1)$ virtually using the information $\{x_{i}(t), \sigma_{i}^{(k)}(t)\}$ sent from the agent $A_{i}$.
The symbol $\numK$ denotes the number of multiple actions.
The $\numK$ virtual inputs $\hat{u}_{i}^{(k)}(t)$ are determined by $\numK$ random variables $\sigma_{i}^{(k)}(t) \in \{-1,1\}^{n}$. 
Each component of the random variable $\sigma_{i}^{(k)}(t)$ independently obeys the Bernoulli distribution with outcome $\pm 1$ equal probabilities.
Let us define $\hat{x}^{(k)}:=[{\hat{x}_{1}^{(k)\mathrm{T}}},..., {\hat{x}_{N}^{(k)\mathrm{T}}}]^{\mathrm{T}}\in \mathbb{R}^{nN}$ and $\sigma^{(k)}(t):=[\sigma_{1}^{(k)}(t)^{\mathrm{T}},...,\sigma_{N}^{(k)}(t)^{\mathrm{T}}]^{\mathrm{T}}$.

Based on the multiple virtual random actions, the global controller $G$ of the PBC law is proposed as follows:
\begin{equation}\label{eq:G_PBC}
\begin{aligned}
\nu(t)&\;=f_{\nu}(x(t),\phi(t),t)
\\&
:=
\begin{bmatrix}
J(\hat{x}^{(1)}(t+1))-J(x(t))\\
\vdots \\
J(\hat{x}^{(\numK)}(t+1))-J(x(t))
\end{bmatrix}
\in  \mathbb{R}^{\numK}
,
\end{aligned}
\end{equation}
where $\nu^{(k)}(t):=J(\hat{x}^{(k)}(t+1))-J(x(t))$ and $\nu(t):=[\nu^{(1)}(t),...,\nu^{(\numK)}(t)]^{\mathrm{T}}$.
Each local controller $L_{i}$ with its state $\phi_{i}(t)$ of the PBC law is proposed as
\begin{equation}\label{eq:L_state_PBC}
\begin{aligned}
&
\phi_{i}(t) :=
\begin{bmatrix}
\phi_{i}^{(1)}(t) \\
\vdots \\
\phi_{i}^{(\numK)}(t) \\
\end{bmatrix}
:=
\begin{bmatrix}
\sigma_{i}^{(1)}(t) \\
\vdots \\
\sigma_{i}^{(\numK)}(t) \\
\end{bmatrix}
\in 
\{-1,1\}^{n\numK}
\end{aligned}
\end{equation}
\begin{equation}\label{eq:L2_PBC}
\begin{aligned}
u_{i}(t)&\;={f}_{u}(\phi_{i}(t),\nu(t),t)
\\&
:=
-a(t)\frac{1}{\numK}\sum_{k=1}^{\numK}\frac{\nu^{(k)}(t)}{c(t)}\phi_{i}^{(k)(-1)}(t)
.
\end{aligned}
\end{equation}
Unlike the BC law, a function ${f}_{\phi}(\phi_{i}(t),\nu(t),t)$ as in (\ref{eq:L1_BC}) is not needed because the state $\phi_{i}(t)$ of the local controller is nothing but the random variables $\sigma_{i}^{(k)}(t)$ $(k=1,...,\numK)$.  
Note that $f_{\nu}$ in (\ref{eq:G_PBC}) is indeed a function of $x(t)$ and $\phi(t)$ because $\hat{x}^{(k)}(t+1)$ is given from $x_{i}(t)$ and $\phi_{i}^{(k)}(t)=\sigma_{i}^{(k)}(t)$, as shown in (\ref{eq:x_predictive}) and (\ref{eq:u_predictive}).
Three important properties of the PBC law are as follows.

First, the PBC is a solution to the main problem and overcomes the limitations of the BC. 
The objective function $J(x)$ is minimized in a manner similar to the BC law.
For the following reasons, the unavailing actions are reduced, and the agents' states are converged at least twice as fast as the BC law under some conditions.
The physical random actions $u_{i}(t)=c(t)\sigma_{i}(t)$ in (\ref{eq:L2_BC}) are no longer used for approximating the gradient of $J(x)$.
The approximation accuracy of the gradient of $J(x)$ is enhanced by taking multiple actions with large $\numK>1$.
In a manner similar to (\ref{eq:E_grad_BC}), the right-hand side of (\ref{eq:L2_PBC}) approximates the gradient as follows:
\begin{equation}\label{eq:E_grad_PBC}
\begin{aligned}
&
\mathrm{E}\Big[ 
\frac{1}{\numK}\sum_{k=1}^{\numK}\frac{\nu^{(k)}(t)}{c(t)}\phi_{i}^{(k)(-1)}(t)
\Big|_{x(t)}
\Big]
=
\partial_{x}J(x(t)) + O(c(t)) 
,
\end{aligned}
\end{equation}
as $c(t) \rightarrow 0$.
The operator $(1/\numK)\sum_{k=1}^{\numK}\{\cdot\}$ in (\ref{eq:E_grad_PBC}) reduces the randomness.
These properties are shown theoretically in Sections \ref{sec_PBC_analysis1} and \ref{sec_PBC_analysis2}.%

Second, the communication volume of the PBC law is slightly greater than that of the BC law.
The global controller $G$ of the PBC law broadcasts $\numK$-dimensional signal $\nu(t)$, whereas the broadcast signal of the BC law is one-dimensional.
In the PBC law, each agent $A_{i}$ transmits the signal $\{x_{i}(t),\phi_{i}(t)\}=\{x_{i}(t), \sigma_{i}^{(1)}(t), ..., \sigma_{i}^{(\numK)}(t)\}$.
The BC law requires that $A_{i}$ transmits only $x_{i}(t)$ or that $G$ observes $x_{i}(t)$ for all $i \in \{1,...,N\}$, in general (see Remark \ref{rem:BC_com_cost}).
However, these differences between the two laws are not critical.
The number $\numK$ of multiple actions is (significantly) smaller than the number $N$ of agents.
The information volume is not so different between $\{x_{i}(t),\phi_{i}(t)\}$ and $x_{i}(t)$ because $\sigma_{i,1}^{(k)}(t) = \pm 1$ is binary ($2$-bit) and is significantly smaller than a quantized real vector $x_{i}(t) \in \mathbb{R}^{n}$.

Third, the advantages of the BC law for large-scale multi-agent systems are preserved in the PBC law.
The PBC law satisfies the three requirements described in Section \ref{sec_target}, which are that communications between agents are not used, that local controllers are identical for all of the agents, and that the broadcast signal is identical.

The following subsections derive the theoretical aspects to show that the PBC law improves the control performance, such as the moving distances of agents and the convergence speed for minimizing an objective function $J(x)$, as compared with the BC law.
Section \ref{sec_PBC_analysis1} analyzes the convergence and performance improvement of the PBC law.
The effect of taking multiple actions with $\numK>1$ is analyzed in Section \ref{sec_PBC_analysis2}.

\subsection{Theoretical analysis of the PBC law: convergence and performance improvement}\label{sec_PBC_analysis1}

This subsection presents a theoretical analysis of the PBC law proposed in Section \ref{sec_PBC_overview}.
In the following, 
Theorem \ref{thm:convergence_of_PBC} shows that the PBC law is a solution to the main problem.
For the case in which the number $\numK$ of multiple random actions is $1$, Theorems \ref{thm:PBC_cost_improvement} and \ref{thm:PBC_dist_reduction} indicate that the control performance of the PBC law is enhanced compared with the BC law.

The condition (c3) is modified for the PBC law:
\begin{enumerate}	
	\item [(c3')]
	The controller gains satisfy $a(t) > 0$ and $c(t) > 0$ for every $t\geq 0$, $\lim_{t \rightarrow \infty}a(t)=0$, $\sum_{t=0}^{\infty}a(t)=\infty$, $\lim_{t \rightarrow \infty}c(t)=0$, $\sum_{t=0}^{\infty}(a(t)/c(t))^{2}<\infty$, and $\sum_{t=0}^{\infty}a(t)c(t)^{2}<\infty$.
\end{enumerate} 
We first show that the PBC law is a solution to the main problem.

\begin{thm}[Convergence of the PBC law]\label{thm:convergence_of_PBC}
For the multi-agent system $\Sigma$, an objective function $J(x)$, and the PBC law in (\ref{eq:G_PBC}), (\ref{eq:L_state_PBC}), and (\ref{eq:L2_PBC}) with any $\numK>0$, if the conditions (c1), (c2), and (c3') hold, then $x(t)$ converges to a (possibly sample-path-dependent) solution to $\partial_{x}J(x) =0$ with probability $1$.
\end{thm}
\begin{pf}
The proof is given in Appendix \ref{pf:convergence_of_PBC}.	
\end{pf}
\begin{rem}
As explained in Remark \ref{rem:BC_solution}, convergence to a solution to $\partial_{x}J(x) =0$ indicates that (\ref{eq:LimJ}) (approximately) holds, and thus the PBC is a solution to the main problem. 
The PBC retains the essential feature of the BC law.
\end{rem}

In the following, we show that the PBC law is superior to the BC law in terms of the convergence speed and the moving distances of agents in the case of $\numK=1$.
The symbols $\{\cdot\}|_{\mathrm{PBC}}$ and $\{\cdot\}|_{\mathrm{BC}}$ denote the results obtained by applying the PBC and BC laws, respectively, where the symbols can be omitted if we discuss either the PBC or the BC.
The following theorem shows that coordination tasks can be performed twice as fast.

\begin{thm}[Achieving tasks at twice the speed] \label{thm:PBC_cost_improvement}
Suppose that $\numK=1$, $a(t)|_{\mathrm{PBC}}=a(2t)|_{\mathrm{BC}}$,  $c(t)|_{\mathrm{PBC}}=c(2t)|_{\mathrm{BC}}$, and $\sigma^{(1)}(t)|_{\mathrm{PBC}}=\sigma(2t)|_{\mathrm{BC}}$ hold for all $t \geq 0$. 
Then, the following relations hold	
	\begin{equation}\label{eq:x_PBC_is_BC}
	\begin{aligned}
	x(t)|_{\mathrm{PBC}}=x(2t)|_{\mathrm{BC}},\quad \forall t \geq 0,
	\end{aligned}
	\end{equation}
	\begin{equation}\label{eq:J_PBC_is_BC}
	\begin{aligned}
	J(x(t))|_{\mathrm{PBC}}=J(x(2t))|_{\mathrm{BC}},\quad \forall t \geq 0.
	\end{aligned}
	\end{equation}
\end{thm}
\begin{pf}
The proof is given in Appendix \ref{pf:PBC_cost_improvement}.	
\end{pf}
\begin{rem}
The speed for achieving the task of the PBC law with $\numK=1$ is twice that of the BC law.
An intuitive reason for this result is that the two-stage transition of the BC law explained in Remark \ref{rem:BC_explanation} is simply combined into a one-stage transition.
\end{rem}

We next analyze the moving distances of agents.
Let us define the total moving distance:
\begin{align}
\totaldist(t)& 
:= \sum_{s=0}^{t-1} \sum_{i=1}^{N}  \| x_{i}(s+1)-x_{i}(s) \|
 = \sum_{s=0}^{t-1} \sum_{i=1}^{N} \| u_{i}(s) \|
.\label{eq:def_dist}
\end{align}
The following theorem shows that the moving distance $\totaldist(t)|_{\mathrm{PBC}}$ of the PBC law is reduced compared with the moving distance $\totaldist(2t)|_{\mathrm{BC}}$ of the BC law in statistical, deterministic, and probabilistic senses.

\begin{thm}[Reduction of the moving distance]\label{thm:PBC_dist_reduction}
Suppose that, $\numK=1$, $a(t)|_{\mathrm{PBC}}=a(2t)|_{\mathrm{BC}}$,  $c(t)|_{\mathrm{PBC}}=c(2t)|_{\mathrm{BC}}$, and $\sigma^{(1)}(t)|_{\mathrm{PBC}}=\sigma(2t)|_{\mathrm{BC}}$ hold for all $t \geq 0$.
Let a binary function $\mathcal{P}(t)$ be $1$ if $J(x)$ is (locally) convex or quasi-convex on the minimum convex set including all possible $x(t+1)|_{\mathrm{BC}}$, and otherwise let be $0$. 
Then, for all $t \geq 1$, the following relations hold:
	\begin{align}
	&
	\mathrm{E}[  \totaldist(2t)|_{\mathrm{BC}} - \totaldist(t)|_{\mathrm{PBC}} ]
	\geq  \sqrt{n} N \sum_{\tau=0}^{t-1} \mathcal{P}(2\tau) c(2\tau)  |_{\mathrm{BC}}
	\nonumber\\&
	\qquad\qquad\qquad\qquad\qquad\;
	\geq 0 
	,\label{eq:E_dist_PBCvsBC}
	\\&
	\totaldist(2t)|_{\mathrm{BC}}  - \totaldist(t)|_{\mathrm{PBC}}
	\geq 0  
	,\;
	\forall \sigma^{(1)}(0)|_{\mathrm{PBC}},..., \sigma^{(1)}(t)|_{\mathrm{PBC}}
	,\label{eq:dist_PBCvsBC}
	\\&
	\mathrm{Pr}\{ \totaldist(2t)|_{\mathrm{BC}} > \totaldist(t)|_{\mathrm{PBC}} \}
	\geq 
	1- \prod_{\tau=0}^{t-1}\Big( 1-\frac{\mathcal{P}(2\tau)}{2} \Big)
	.\label{eq:dist_PBCvsBC_prob}
	\end{align}
\end{thm}
\begin{pf}
The proof is given in Appendix \ref{pf:PBC_dist_reduction}.
\end{pf}
\begin{rem}
The local convexity of $J(x)$ increases the effectiveness of the PBC law because the binary function $\mathcal{P}(2t)$ is activated from $0$ to $1$ by the local convexity.	
If $J(x)$ is convex or quasi-convex for all $x$, (\ref{eq:E_dist_PBCvsBC}) and (\ref{eq:dist_PBCvsBC_prob}) reduce to
$\mathrm{E}[  \totaldist(2t)|_{\mathrm{BC}} - \totaldist(t)|_{\mathrm{PBC}} ]
\geq   \sqrt{n} N \sum_{\tau=0}^{t-1} c(2\tau)  |_{\mathrm{BC}}$ and $\mathrm{Pr}\{ \totaldist(2t)|_{\mathrm{BC}} > \totaldist(t)|_{\mathrm{PBC}} \} \geq 1-(1/2^{t})$, respectively.
\end{rem}
\begin{rem}
In Theorem \ref{thm:PBC_dist_reduction}, it is reasonable that the time scales of the PBC and BC laws are normalized as $t$ and $2t$, respectively, because $x(t)|_{\mathrm{PBC}}=x(2t)|_{\mathrm{BC}}$ and $J(x(t))|_{\mathrm{PBC}}=J(x(2t))|_{\mathrm{BC}}$ hold in Theorem \ref{thm:PBC_cost_improvement}.  
\end{rem}

\subsection{Theoretical analysis of the PBC law: performance improvement by multiple actions}\label{sec_PBC_analysis2}

This subsection presents a theoretical analysis of the PBC law for the case of taking multiple actions, i.e., $\numK>1$.
Theorems \ref{thm:MPBC_per_step_improvement} and \ref{thm:MPBC_multi_step_improvement} and Proposition \ref{thm:quad_J} give results for control performance improvement with $\numK>1$.

To the best of our knowledge, few studies with respect to SPSA using multiple random numbers have been analyzed theoretically.
The asymptotic properties of SPSA using multiple random numbers are discussed in \cite{Spall92}.
It is also well known that the variance of the average of $\numK$ samples is decreased $1/\numK$ times compared to the variance of one sample.
This property is applied to the approximate gradient of $J(x)$, i.e., the input in (\ref{eq:L2_PBC}):
\begin{equation}
\begin{aligned}
\mathrm{Cov} [  u_{i}(t)|_{x(t)} ]
=
\frac{1}{\numK}
\big(
\mathrm{Cov} [  u_{i}(t)|_{x(t)} ] 
\big|_{\numK=1}
\big)
.
\end{aligned}
\end{equation}
Using a large $\numK$ enhances the approximation accuracy of the gradient in terms of the covariance.

In the following, some theorems are derived to discuss the influence of $\numK>1$.
We derive a theorem to show that increasing the value of $\numK$ is efficient for enhancing the control performance per step of the PBC law.

\begin{thm}[Performance improvement per step]\label{thm:MPBC_per_step_improvement}	
For a given $x(t)$, if the PBC law is applied, the following relations hold 
\begin{enumerate}
\item 
If $J(x)$ is (locally) convex on the minimum convex set including all possible $x(t+1)$ for the given $x(t)$,
the reduction of $J(x)$ per step is enhanced by increasing $\numK$
\begin{equation}\label{eq:per_step_cost_imporvvement_by_MPBC}
\begin{aligned}
&
\forall 1 \leq \numK_{b} < \numK_{a} < \infty
,
\\&
\mathrm{E}[J(x(t+1)) |_{x(t)} ] |_{\numK=\numK_{a}} 
{\leq}
\mathrm{E}[J(x(t+1)) |_{x(t)} ] |_{\numK=\numK_{b}} 
,
\end{aligned}
\end{equation}	
where the equality in (\ref{eq:per_step_cost_imporvvement_by_MPBC}) is possible only if $J(x)$ is not strictly convex but convex.
The inequality sign in (\ref{eq:per_step_cost_imporvvement_by_MPBC}) is reversed if $J(x)$ is concave.

\item For any $\kappa \geq 1 \in \mathbb{R}$, the power of the moving distance per step is reduced by increasing $\numK$
\begin{equation}\label{eq:per_step_distance_reduce_by_MPBC}
\begin{aligned}
&
\forall 1 \leq \numK_{b} < \numK_{a} < \infty
,
\\&
\mathrm{E}\Big[  \sum_{i=1}^{N} \|x_{i}(t+1)-x_{i}(t)\|^{\kappa}  \Big|_{x(t)} \Big] \Big|_{\numK=\numK_{a}} 
\\&
\leq
\mathrm{E}\Big[  \sum_{i=1}^{N} \|x_{i}(t+1)-x_{i}(t)\|^{\kappa}  \Big|_{x(t)} \Big] \Big|_{\numK=\numK_{b}} 
,
\end{aligned}
\end{equation}	
where the equality in (\ref{eq:per_step_distance_reduce_by_MPBC}) is possible only if $\kappa=1$.

\end{enumerate}
\end{thm}
\begin{pf}	
The proof is given in Appendix \ref{pf:MPBC_per_step_improvement}.	
\end{pf}
\begin{rem}
This theorem means that the per-step control performance in terms of the speed for minimizing $J(x)$ and the moving distance is statistically improved by using a large $\numK$.
If $J(x)$ is locally concave around $x(t)$, setting $\numK=1$ is suitable for minimizing $J(x)$ quickly.
\end{rem}

Theorem \ref{thm:MPBC_per_step_improvement} focuses on the control performance at each step, which is independent of other steps.
Based on Theorem \ref{thm:MPBC_per_step_improvement}, we derive a result for a suitable objective function $J(x)$ such that the control performance over multiple steps is improved.

\begin{thm}[Performance over multiple steps]\label{thm:MPBC_multi_step_improvement}	
For any $t \geq 1$ and a given $x(0)$, let $\mathbb{S}_{x}$ be the minimum convex set including all possible $x(0)$, $x(1)$,..., and $x(t)$.
Suppose that, for all $\sigma^{(k)}(m)$ $(m \in \{0,...,t-1\},\; k \in \{1,...,\numK\})$ and for all $s \in \{0,...,t\}$, $J_{s,t}(\tilde{x}):=J(x(t))|_{x(s)=\tilde{x}}$ and $\totaldist_{s,t}(\tilde{x}):=(\totaldist(t)-\totaldist(s))|_{x(s)=\tilde{x}}$ are (locally) convex in $\tilde{x} \in \mathbb{S}_{x}$.	
The PBC law enhances the reductions of $J(x(t))$ and  $\totaldist(t)$ by increasing $\numK$
\begin{align}
&
\forall 1 \leq \numK_{b} < \numK_{a} < \infty,
\nonumber\\& \quad
\mathrm{E}[ J(x(t))  ] |_{\numK=\numK_{a}} 
\leq
\mathrm{E}[ J(x(t))  ] |_{\numK=\numK_{b}} 
,\label{eq:multisteps_cost_imporvvement_by_MPBC}
\\&
\forall 1 \leq \numK_{b} < \numK_{a} < \infty,
\nonumber\\& \quad
\mathrm{E}[ \totaldist(t) ] |_{\numK=\numK_{a}} 
\leq
\mathrm{E}[ \totaldist(t) ] |_{\numK=\numK_{b}} 
,\label{eq:multisteps_distance_imporvvement_by_MPBC}
\end{align}
where the equalities in (\ref{eq:multisteps_cost_imporvvement_by_MPBC}) and (\ref{eq:multisteps_distance_imporvvement_by_MPBC}) are possible only if $J_{s,t}(\tilde{x})$ and $\totaldist_{s,t}(\tilde{x})$ are not strictly convex but convex in $\tilde{x}$, respectively.
\end{thm}
\begin{pf}
The proof is given in Appendix \ref{pf:MPBC_multi_step_improvement}.
\end{pf}
\begin{rem}
Although Theorem \ref{thm:MPBC_multi_step_improvement} considers the time steps from $0$ to $t$, it can be applied to the time steps from $\tau$ to $\tau+t$  for a given $x(\tau)$ and $\tau \geq 0$ in a similar manner.
\end{rem}

Theorem \ref{thm:MPBC_multi_step_improvement} shows that the performance of $J(x(t))$ and $\totaldist(t)$ are enhanced over multiple steps for suitable objective functions $J(x)$.
There are some choices of $J(x)$ satisfying the assumptions in Theorem \ref{thm:MPBC_multi_step_improvement}. 
One such class of $J(x)$ is the convex quadratic function class, as described below.

\begin{prop}[Convex quadratic functions]\label{thm:quad_J}
If $J(x)$ is convex and quadratic in $x$, the assumptions in Theorem \ref{thm:MPBC_multi_step_improvement} hold for any $x(0)$ (it is possible that $J_{s,t}(\tilde{x})$ and $\totaldist_{s,t}(\tilde{x})$ is not strictly convex but convex).
\end{prop}
\begin{pf}
The proof is given in Appendix \ref{pf:quad_J}.	
\end{pf}
\begin{rem}
	Although Proposition \ref{thm:quad_J} focuses on convex quadratic functions, the PBC law can be successfully applied to various non-convex objective functions, which are locally convex and/or quadratic. 
	Section \ref{sec_simulation} demonstrate the PBC law for non-convex objective functions.
\end{rem}

Theorem \ref{thm:MPBC_multi_step_improvement} and Proposition \ref{thm:quad_J} indicate that convex quadratic objective functions are more suitable than other objective functions. 
Some locally/globally convex quadratic objective functions for coordination tasks are introduced in the following.

\begin{enumerate}
\item 
\textbf{Coverage control.}
The following objective function for coverage control is locally convex and quadratic in $x$ 
\begin{align} \label{eq:obj_coverage}
&
J_{\mathrm{obj}}(x)=\frac{\mathbb{V}_{q}}{N_{q}}\sum_{j=1}^{N_{q}} \min_{i\in \{1,...,N\}} \| q_{j} - x_{i}\|^{2}
,
\end{align}
where $q_{j}$ are uniformly distributed on the state space, and $\mathbb{V}_{q}:=\int_{\mathbb{S}_{q}} 1 \mathrm{d}q$ is the $n$-dimensional volume of a workspace $\mathbb{S}_{q} \subset \mathbb{R}^{n}$ $(q \in \mathbb{S}_{q})$. 
The operator $({\mathbb{V}_{x}}/{N_{q}}) \sum_{j=1}^{N_{q}}\{ \cdot \}$ approximates $\int_{\mathbb{S}_{q}} \{ \cdot \} \mathrm{d}q$.

\item 
\textbf{Rendezvous control with formation selection.}
The following objective function for a rendezvous control is locally convex and quadratic in $x$
\begin{align} \label{eq:obj_rendezvous}
&
J_{\mathrm{obj}}(x)=\min_{\theta \in \mathbb{S}_{\theta} }  \frac{1}{N^{2}}   \sum_{i=1}^{N}\sum_{j=1}^{N} \|x_{i}-x_{j} - r_{i,j}(\theta) \|^{2}
,
\end{align}
because, for all $i$ and $j$, the L2-norms are quadratic convex and their sum is quadratic convex.
The symbol $r_{i,j}(\theta) \in \mathbb{R}^{n}$ is a relative target position between the agents $i$ and $j$.
Various target formations are constructed from $r_{i,j}(\theta)$ parameterized by a global formation parameter $\theta \in \mathbb{S}_{\theta}$.
Assuming that $\mathbb{S}_{\theta}$ is finite, an optimal $\theta$ minimizing $J_{\mathrm{obj}}(x)$ is selected from among $\mathbb{S}_{\theta}$. 
The PBC law can select such an (local) optimal formation with $\theta$ by minimizing $J_{\mathrm{obj}}(x)$, whereas standard distributed control cannot optimize $\theta$.
If $\mathbb{S}_{\theta}$ includes only one element, $J_{\mathrm{obj}}(x)$ in (\ref{eq:obj_rendezvous}) reduces to  a global convex objective function for well-known rendezvous control.

\item 
\textbf{Assignment control.}
For given target locations $y_{i} \in \mathbb{R}^{n}$ $(i=1,...,N)$, the following objective function for an assignment task is globally convex and quadratic in $x$
\begin{align} \label{eq:obj_assignment}
&
J_{\mathrm{obj}}(x,\AssignID)=\sum_{i=1}^{N} \|x_{i}-y_{\AssignID_{i}}\|^{2} 
,
\end{align}
where $\AssignID_{i} \in \{1,...,N\}$ and $\AssignID_{i} \neq \AssignID_{j}$ for any $i \neq j$.
The indices $\AssignID:=[\AssignID_{1},...,\AssignID_{N}]$ indicate the pairs of target locations and agents to be assigned. 
Some assignment algorithms determine suitable indices.
For example, the Hungarian algorithm gives optimal indices satisfying $\argmin_{\AssignID}J_{\mathrm{obj}}(x,\AssignID)$ with the complexity $O(N^{3})$ \cite{Munkres57}.

\end{enumerate}
The above-mentioned functions in (\ref{eq:obj_coverage}) and (\ref{eq:obj_rendezvous}) include the operator $\min\{\cdot\}$, which invokes indifferentiable points on the state space.
If the $C^{2}$ continuity in the condition (c1) is strictly concerned, a modified log-sum-exp form\footnote{%
	The log-sum-exp form approximates the max function: $\max_{j} \SmallParamCtwo f_{j} \leq  \ln \sum_{j=1}^{n_{j}} \exp  \SmallParamCtwo f_{j}  \leq \max_{j} \SmallParamCtwo f_{j} + \ln n_{j}$ (see Section 3.1.5 in \cite{Boyd04}). 
	For $\SmallParamCtwo<0$, this inequality is transformed as $ \SmallParamCtwo \min_{j} f_{j} \leq \ln \sum_{j=1}^{n_{j}}  \exp \SmallParamCtwo f_{j}  \leq \SmallParamCtwo \min_{j} f_{j} + \ln n_{j}$, which corresponds to (\ref{eq:log_sum_exp}). 
}%
can approximate the operator $\min_{j}\{\cdot\}$ by a smooth function $(1/\SmallParamCtwo) \ln \sum_{j} \exp \SmallParamCtwo \{ \cdot \}$
\begin{align}
\frac{1}{\SmallParamCtwo}  \ln n_{j}
&\leq
\frac{1}{\SmallParamCtwo} \ln  \sum_{j=1}^{n_{j}} \exp \SmallParamCtwo f_{j} 
-\min_{ j \in \{1,...,n_{j}\} } f_{j} 
\leq
0
,\label{eq:log_sum_exp}
\end{align}
for any $f_{j} \in \mathbb{R}$ $(j=1,2,...,n_{j})$, where $\SmallParamCtwo <0 \in \mathbb{R}$.
The term $(1/\SmallParamCtwo) \ln n_{j} \leq 0$ is the approximation residual.
Sufficiently small $\SmallParamCtwo$ realizes a better approximation, i.e., $(1/\SmallParamCtwo) \ln n_{j} \rightarrow 0$ as $\SmallParamCtwo \rightarrow - \infty$.

\section{Numerical example}\label{sec_simulation}

\definecolor{myPlotB}{rgb}{ 0,0,1 }
\definecolor{myPlotG}{rgb}{ 0,0.7,0 }
\definecolor{myPlotR}{rgb}{ 1,0,0 }
\definecolor{myPlotY}{rgb}{ 0.7,0.7,0 }

In this section, the effectiveness of the PBC law is evaluated in comparison with the BC law.
The two types of multi-agent coordination tasks are demonstrated.

\subsection{Settings}\label{sec_sim_setting}

\newcommand{\terminaltime}{300}

The number of agents is set to $N=15$ in the two-dimensional state space, i.e., $n=2$.
The parameters in (\ref{eq:barrier_switch}) were respectively set to $l_{1} = 100$ and $l_{2} = 101$, which determine the size of the workspace.
The terminal time of the simulation is set to $\terminaltime$.
The controller gains are set such that the conditions (c3) and (c3') are satisfied
\begin{align}
	a(t)|_{\mathrm{PBC}}
	&=\frac{a_{0}}{(t +  t_{v} )^{a_p}} 
	, \label{eq:def_gain_a} \\
	c(t)|_{\mathrm{PBC}}
	&=\frac{c_{0}}{(t +  t_{v} )^{c_p}} 
	, \label{eq:def_gain_c}
\end{align}
where $a(2t)|_{\mathrm{BC}}=a(2t+1)|_{\mathrm{BC}}=a(t)|_{\mathrm{PBC}}$ and $c(2t)|_{\mathrm{BC}}=c(2t+1)|_{\mathrm{BC}}=c(t)|_{\mathrm{PBC}}$ in the case of the BC law.
Based on conditions for the convergence/divergence of the $p$-series (see Section 8.1.2 \cite{Polyanin06}), the conditions (c3) and (c3') are satisfied if $t_{v}>0$, $0 < a_{p}\leq 1$, $c_{p}>0$, $2a_{p} - 2c_{p} > 1$, and $a_{p} + 2 c_{p} > 1$ hold.
According to these conditions, the coefficients were set to $a_{0}=2$, $a_{p}=0.7$, $c_{0}=0.003$, $c_{p}=0.16$, and $t_{v}=20$.

\subsection{Broadcast rendezvous control with formation selection \label{sec_sim_rendezvous}}

\begin{figure}[t]	
	\begin{minipage}[b]{.5\linewidth}
		\centering
		\includegraphics[width=0.95\linewidth]{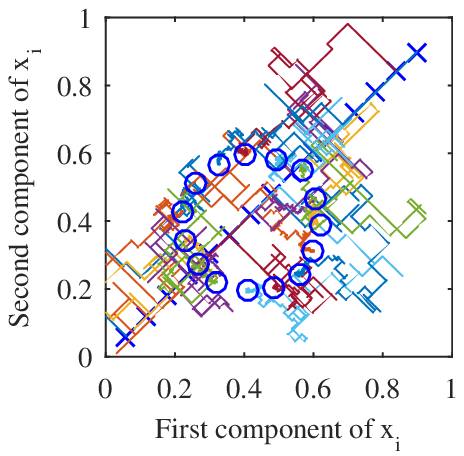}
		\vspace*{-0.082in}
		\subcaption{BC.}  
	\end{minipage}%
	\begin{minipage}[b]{.5\linewidth}
		\centering
		\includegraphics[width=0.95\linewidth]{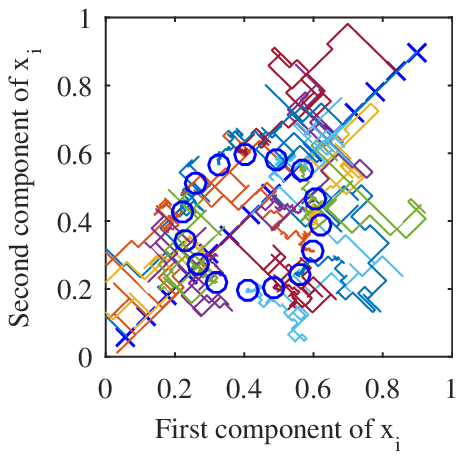}
		\vspace*{-0.082in}
		\subcaption{PBC ($\numK=1$).}  
	\end{minipage}%
	
	\vspace*{+0.051in}
	
	\begin{minipage}[b]{.5\linewidth}
		\centering
		\includegraphics[width=0.95\linewidth]{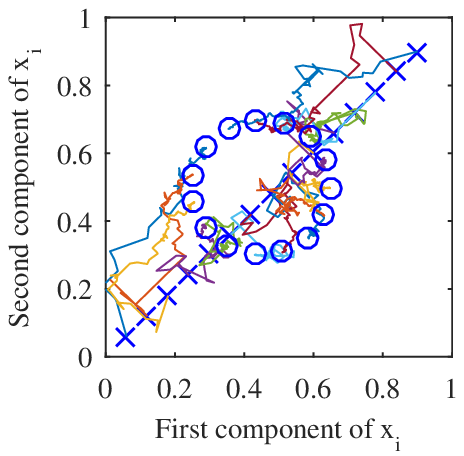}
		\vspace*{-0.082in}
		\subcaption{PBC ($\numK=3$).}  
	\end{minipage}%
	\begin{minipage}[b]{.5\linewidth}
		\centering
		\includegraphics[width=0.95\linewidth]{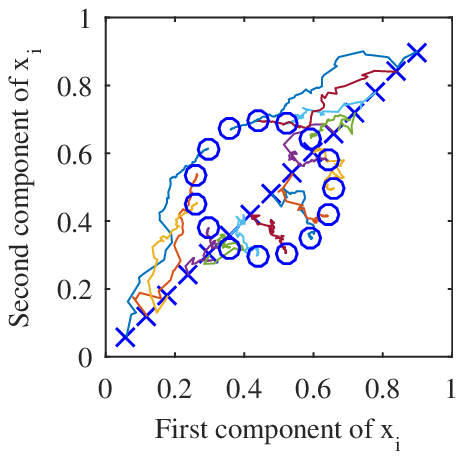} 
		\vspace*{-0.082in}
		\subcaption{PBC ($\numK=10$).} 
	\end{minipage}%
	\caption{State transitions at the terminal time $t=\terminaltime$ for the rendezvous control}
	\label{fig:state_rendezvous}
\end{figure}

\begin{figure}[t]	
	{\scriptsize
		\vspace*{0.08in}
		\hfill
		\tabcolsep = 1mm
		\renewcommand\arraystretch{1.2}
		\begin{tabular}{|l l l  l|}
			\hline
			\textcolor{black}{$\blacksquare$}:BC &
			\textcolor{myPlotB}{$\blacksquare$}:PBC ($\numK=1$) &
			\textcolor{myPlotG}{$\blacksquare$}:PBC ($\numK=3$) &
			\textcolor{myPlotR}{$\blacksquare$}:PBC ($\numK=10$) 
			\\\hline				
		\end{tabular}	
	}

	\vspace*{+0.051in}

	\begin{minipage}[b]{.5\linewidth}
		\centering
		\includegraphics[width=0.95\linewidth]{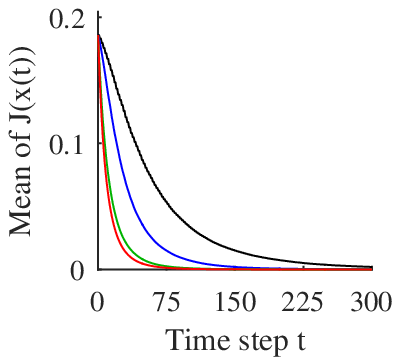} 
		\vspace*{-0.082in}
		\subcaption{Mean of $J(x(t))$.} 
	\end{minipage}%
	\begin{minipage}[b]{.5\linewidth}
		\centering
		\includegraphics[width=0.95\linewidth]{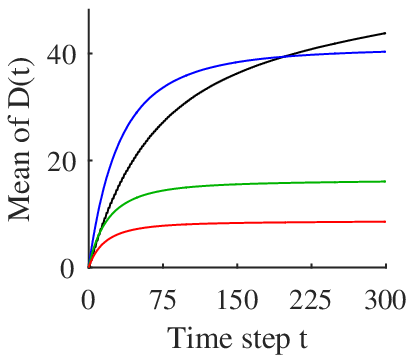}
		\vspace*{-0.082in}
		\subcaption{Mean of $\totaldist(t)$.}  
	\end{minipage}%
	
	\vspace*{+0.051in}
	
	\begin{minipage}[b]{.5\linewidth}
		\centering
		\includegraphics[width=0.95\linewidth]{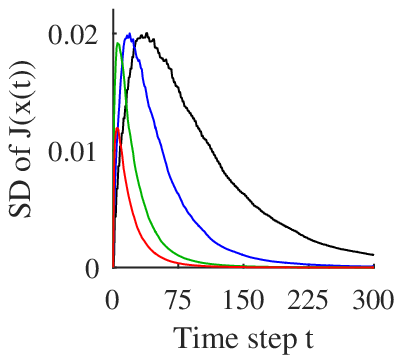}
		\vspace*{-0.082in}
		\subcaption{SD of $J(x(t))$.}  
	\end{minipage}%
	\begin{minipage}[b]{.5\linewidth}
		\centering
		\includegraphics[width=0.95\linewidth]{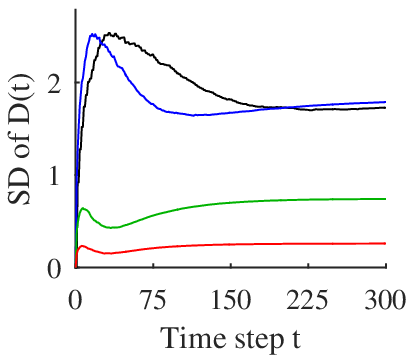}
		\vspace*{-0.082in}
		\subcaption{SD of $\totaldist(t)$.}  
	\end{minipage}%
	\caption{Transitions of the objective function $J(x(t))$ and the moving distance $\totaldist(t)$ for the rendezvous control}
	\label{fig:results_J_D_rendezvous}
\end{figure}

This subsection evaluates a rendezvous control task with the objective function $J_{\mathrm{obj}}(x)$ in (\ref{eq:obj_rendezvous}).
Recall that minimizing $J_{\mathrm{obj}}(x)$ in (\ref{eq:obj_rendezvous}) automatically selects a best formation parameterized by $\theta \in \mathbb{S}_{\theta}$.
The state trajectories are shown in Fig. \ref{fig:state_rendezvous}.
The $\times$ and $\circ$ symbols denote the initial states $x_{i}(0)$ and the terminal states $x_{i}(\terminaltime)$, respectively. 
The colored lines represent the trajectories of the agents.
As illustrated in Fig. \ref{fig:state_rendezvous}, the initial state was set to $x_{i}(0)=(0.9 i/N)[1,1]^{\mathrm{T}}$.
The target formations were set such that $r_{i,j}(\theta)=y_{i}(\theta)-y_{j}(\theta)$ in (\ref{eq:obj_rendezvous}), where $y_{i}(\theta):=0.2 [\cos (2 \pi (i + \theta) / N) ,  \sin (2 \pi (i + \theta) / N)]^{\mathrm{T}}$ and $\theta \in \mathbb{S}_{\theta} := \{1,...,N\}$. 
In Fig. \ref{fig:state_rendezvous}, the BC law and the PBC law with $\numK=1$ caused the large random actions for the agents.
Applying large $\numK \in \{3,10\}$ yields smooth trajectories, reducing the unavailing random actions.
Figure \ref{fig:results_J_D_rendezvous} shows the transitions of the objective function $J(x(t))$ and the moving distance $\totaldist(t)$ with respect to the mean and the standard deviation (SD) for $500$ trials (with different random seeds).
The both $J(x(t))$ and $\totaldist(t)$ were successfully reduced by the PBC law.
Their means and SDs were decreased more by using large $\numK \in \{3,10\}$ as compared with those of the BC law and the PBC law with $\numK=1$.

\subsection{Broadcast coverage control \label{sec_sim_coverage}}

\begin{figure}[t]	
	\begin{minipage}[b]{.5\linewidth}
		\centering
		\includegraphics[width=0.95\linewidth]{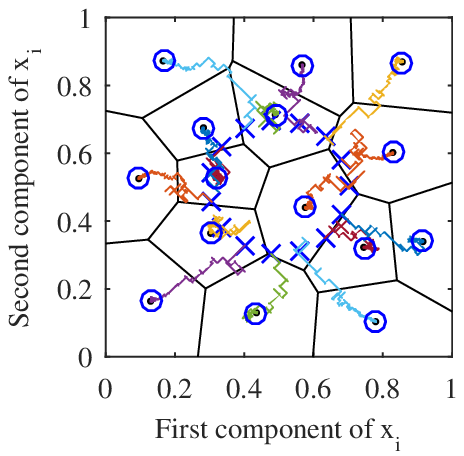}
		\vspace*{-0.082in}
		\subcaption{BC.}  
	\end{minipage}%
	\begin{minipage}[b]{.5\linewidth}
		\centering
		\includegraphics[width=0.95\linewidth]{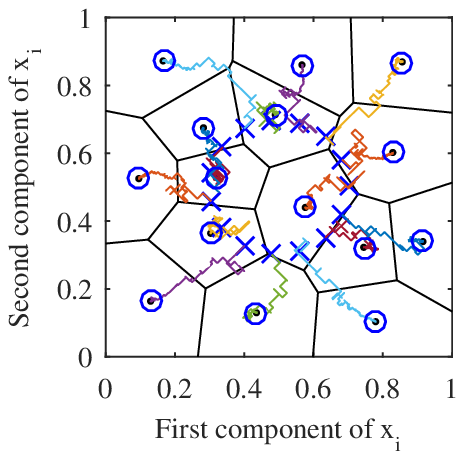}
		\vspace*{-0.082in}
		\subcaption{PBC ($\numK=1$).}  
	\end{minipage}%
	
	\vspace*{+0.051in}
	
	\begin{minipage}[b]{.5\linewidth}
		\centering
		\includegraphics[width=0.95\linewidth]{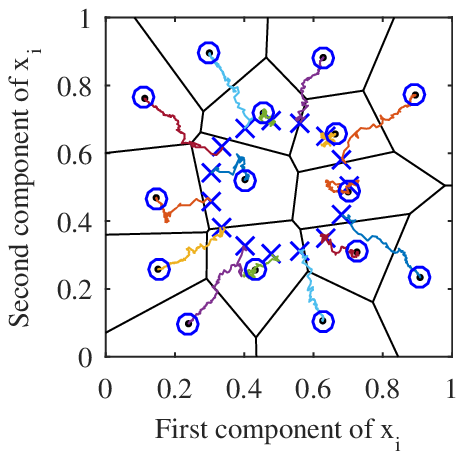}
		\vspace*{-0.082in}
		\subcaption{PBC ($\numK=3$).}  
	\end{minipage}%
	\begin{minipage}[b]{.5\linewidth}
		\centering
		\includegraphics[width=0.95\linewidth]{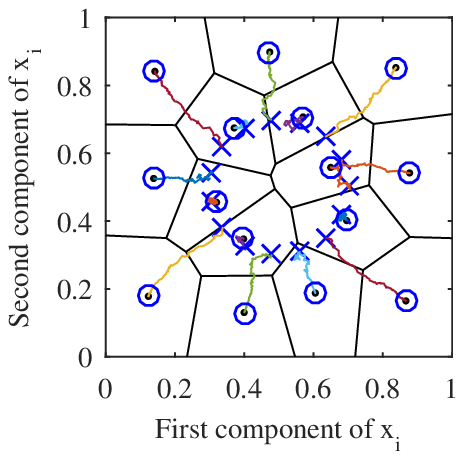}
		\vspace*{-0.082in}
		\subcaption{PBC ($\numK=10$).}  
	\end{minipage}%
	\caption{State transitions at the terminal time $t=\terminaltime$ for the coverage control}
	\label{fig:state_coverage}
\end{figure}

\begin{figure}[t]	
	{\scriptsize
		\vspace*{0.08in}
		\hfill
		\tabcolsep = 1mm
		\renewcommand\arraystretch{1.2}
		\begin{tabular}{|l l l  l|}
			\hline
			\textcolor{black}{$\blacksquare$}:BC &
			\textcolor{myPlotB}{$\blacksquare$}:PBC ($\numK=1$) &
			\textcolor{myPlotG}{$\blacksquare$}:PBC ($\numK=3$) &
			\textcolor{myPlotR}{$\blacksquare$}:PBC ($\numK=10$) 
			\\\hline				
		\end{tabular}	
	}	
	
	\vspace*{+0.051in}
	
	\begin{minipage}[b]{.5\linewidth}
		\centering
		\includegraphics[width=0.95\linewidth]{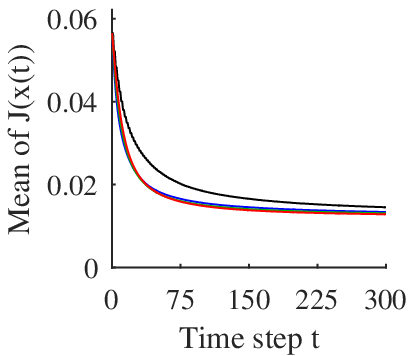}
		\vspace*{-0.082in}
		\subcaption{Mean of $J(x(t))$.}  \label{fig:mean_J_coverage} 
	\end{minipage}%
	\begin{minipage}[b]{.5\linewidth}
		\centering
		\includegraphics[width=0.95\linewidth]{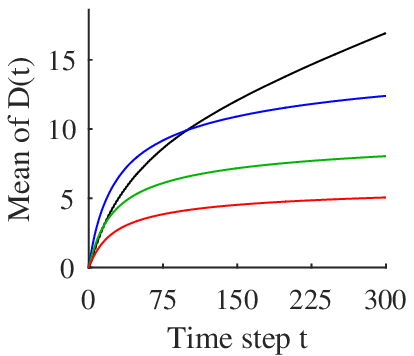}
		\vspace*{-0.082in}
		\subcaption{Mean of $\totaldist(t)$.}  
	\end{minipage}%
	
	\vspace*{+0.051in}
	
	\begin{minipage}[b]{.5\linewidth}
		\centering
		\includegraphics[width=0.95\linewidth]{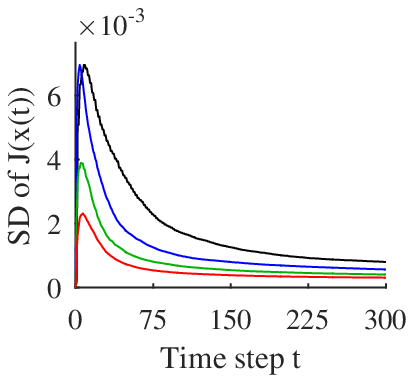}
		\vspace*{-0.082in}
		\subcaption{SD of $J(x(t))$.}  
	\end{minipage}%
	\begin{minipage}[b]{.5\linewidth}
		\centering
		\includegraphics[width=0.95\linewidth]{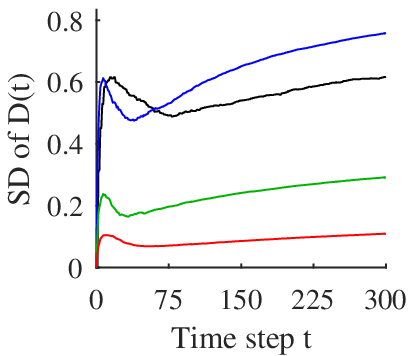}
		\vspace*{-0.082in}
		\subcaption{SD of $\totaldist(t)$.}  
	\end{minipage}%
	\caption{Transitions of the objective function $J(x(t))$ and the moving distance $\totaldist(t)$ for the coverage control}
	\label{fig:results_J_D_coverage}
\end{figure}

A coverage control task is evaluated with the objective function $J_{\mathrm{obj}}(x)$ in (\ref{eq:obj_coverage}), where $q_{j}$ are sampled on the space $[0,1] \times [0,1]$ at $0.01$ intervals in each direction.
Figure \ref{fig:state_coverage} shows the state trajectories of the PBC law with various $\numK \in \{1,3,10\}$ compared with the BC law.
The initial state was set to $x_{i}(0)=[0.5 + 0.2 \cos (2 \pi i / N), 0.5 + 0.2 \sin (2 \pi i / N)]^{\mathrm{T}}$, as illustrated in Fig. \ref{fig:state_coverage}.
The black lines indicate the Voronoi diagram composed by the agents.
We see that increasing $\numK$ reduces the unavailing random actions of the agents.
Figure \ref{fig:results_J_D_coverage} shows $J(x(t))$ and $\totaldist(t)$ in terms of the mean and the SD for $500$ trials (with different random seeds).
In Fig. \ref{fig:results_J_D_coverage}\subref{fig:mean_J_coverage}, $J(x(t))$ was almost the same regardless of the value of $\numK$.
It appears that $J(x(t))$ is not convex around the initial state $x(0)$.
Recall Theorem \ref{thm:MPBC_per_step_improvement}, which indicates that increasing $\numK$ does not ensure the improvement of $J(x(t))$ for non-convex $J(x(t))$. 
Nevertheless, the PBC law reduces both $J(x(t))$ and $\totaldist(t)$ compared with the existing BC law.

\section{Conclusion}\label{sec_conclusion}

In the present paper, we proposed the PBC law for achieving multi-agent coordination tasks with low communication volume without any agent-to-agent communication.
To address drawbacks of the BC law discussed in Section \ref{sec_main_problem}, a simple but efficient solution was proposed in Section \ref{sec_PBC_overview}.
The solution is to take multiple virtual actions of agents rather than a single physical random action.
Sections \ref{sec_PBC_analysis1} and \ref{sec_PBC_analysis2} proved the theoretical aspects of the PBC law.
Theorem \ref{thm:convergence_of_PBC} proved asymptotic convergence of multi-agent coordination with probability $1$.
Theorems \ref{thm:PBC_cost_improvement} and \ref{thm:PBC_dist_reduction} show that the PBC law is superior to the BC law in terms of the convergence speed and the moving distance in the case of taking a single action ($\numK=1$).
In the case of multiple actions ($\numK>1$), Theorem \ref{thm:MPBC_per_step_improvement} showed that the performance improvements per step become more effective by increasing the number $\numK$ of multiple actions.
Theorem \ref{thm:MPBC_multi_step_improvement} proved that the performance improvements are retained over multiple steps for suitable objective functions.
Such functions are convex quadratic in the state shown in Proposition \ref{thm:quad_J}.
Section \ref{sec_simulation} demonstrated that the PBC law is superior to the BC law through two types of coordination tasks.

In the future, we intend to extend the PBC law for nonlinear and/or stochastic multi-agent systems.
The PBC law will be applied to merging tasks of multiple vehicles on congested roads.


\appendix

\section{Proof of Theorem \ref{thm:convergence_of_PBC}} \label{pf:convergence_of_PBC}
This theorem is proven in a manner similar to the case of the BC law \cite{azuma2013broadcast}.
Unfortunately, the difficulty in the case of the PBC law arises from taking multiple trials ($\numK>1$).
The proof of the convergence of the BC law in Appendices A.1 and A.2 of \cite{azuma2013broadcast} is summarized (with slight modification) as follows.
\begin{lem}[Stochastic systems for the BC \cite{azuma2013broadcast}]\label{thm:convergence_proof_of_BC}
	Let us define $t_{s}:=2s$ for $s \in \{0,1,...\}$. 
	Supposing that the conditions (c1), (c2), and (c3) hold, $x(t_{s})$ converges to a (possibly sample-path-dependent) solution to $\partial_{x}J(x) =0$ with probability $1$ if the following conditions hold:
	\begin{enumerate}
		\item [(c4)] For all $t_{s}$, $x(t_{s})$ obeys
		\begin{equation}\label{eq:sys_for_Borkar_lem_BC_ver}
		\begin{aligned}
		&
		x(t_{s+1})
		= x(t_{s})
		\\&
		 - a(t_{s}) \{ \partial_{x}J(x(t_{s})) + \BCnoiseRand(t_{s+1}) + \BCnoiseOffset(t_{s+1}) \}
		,
		\end{aligned}
		\end{equation}
		where $\BCnoiseRand(t_{s}) \in \mathbb{R}^{nN}$ is a random vector and 
		$\BCnoiseOffset(t_{s}) \in \mathbb{R}^{nN}$ is a bounded random vector satisfying $\BCnoiseOffset(t_{s}) \rightarrow 0$ as $t_{s} \rightarrow \infty$.
		
		\item [(c5)]	
		$\sup_{t_{s}} \|x(t_{s})\| < \infty$ holds with probability 1. 
		
		\item [(c6)]	
		For any $m \in \{0,1,...\}$, the stochastic process $\{ \sum_{s=0}^{m} a(t_{s})\BCnoiseRand(t_{s+1}) \}$ is a square integrable martingale, and $\sum_{s=0}^{m} a(t_{s})^{2} \mathrm{E}[ \|\BCnoiseRand(t_{s+1})\|^{2} | \mathcal{F}_{s}]<\infty$ holds with probability $1$ for the filtration $\mathcal{F}_{s}$ generated by 
		$\{x(t_{0}), x(t_{s^{\prime}}),\BCnoiseRand(t_{s^{\prime}}),\BCnoiseOffset(t_{s^{\prime}}) | 1 \leq s^{\prime} \leq s \}$.

	\end{enumerate}
\end{lem}	

We apply Lemma \ref{thm:convergence_proof_of_BC} to the PBC law.
Let us define the following functions for brief notation:
\begin{align}
&
\defGrad(t,\sigma)
:=
\frac{J(x(t)+c(t)\sigma) - J(x(t))}{c(t)}\sigma^{(-1)}
, \label{eq:input_PBC}
\\&
\ERRdefGrad(t)
:=
\frac{1}{\numK}\sum_{k=1}^{\numK} \defGrad(t,\sigma^{(k)}(t))
-
\mathrm{E}\Big[
\frac{1}{\numK}\sum_{k=1}^{\numK} \defGrad(t,\sigma^{(k)}(t))
\Big|_{x(t)}
\Big]
.\label{eq:def_error_of_grad}
\end{align}
By substituting (\ref{eq:G_PBC}), (\ref{eq:L_state_PBC}), (\ref{eq:L2_PBC}), (\ref{eq:E_grad_PBC}), and (\ref{eq:def_error_of_grad}) into (\ref{eq:Agent}), the transition of $x(t)$ under the PBC law is given by
\begin{equation}\label{eq:transition_PBC}
\begin{aligned}
x(t+1)&=x(t)- a(t)\frac{1}{\numK}\sum_{k=1}^{\numK} \defGrad(t,\sigma^{(k)}(t))
\\&
=
x(t) - a(t) \{ \partial_{x}J(x(t)) + \ERRdefGrad(t) + O(c(t)) \}
,
\end{aligned}
\end{equation}
as $c(t) \rightarrow 0$.
If the condition (c5) in  Lemma \ref{thm:convergence_proof_of_BC} holds, then the term $O(c(t))$ is bounded and $O(c(t))\rightarrow 0$ as $t \rightarrow \infty$ from the condition (c3').
The pair of the condition (c3) and the definition $t_{s}:=2s$ in Lemma \ref{thm:convergence_proof_of_BC}  can be replaced by the pair of  the condition (c3') and $t_{s}:=s$.
The condition (c4) then holds under the condition (c5) by regarding $t$, $t+1$, $\ERRdefGrad(t)$, and $O(c(t))$ as $t_{s}$, $t_{s+1}$, $\BCnoiseRand(t_{s+1})$, and $\BCnoiseOffset(t_{s+1})$.
We next prove the condition (c5).
The system dynamics in (\ref{eq:transition_PBC}) is transformed as follows:
\begin{align}
\|x(t+1)\|^{2}&=\Big\| \frac{1}{\numK} \sum_{k=1}^{\numK} \{ x(t)- a(t) \defGrad(t,\sigma^{(k)}(t)) \} \Big\|^{2}
\nonumber \\&
\leq
\| x(t)- a(t) \defGrad(t,\sigma^{(k_{*}(t))}(t))  \|^{2}
,\label{eq:bounded_state_transition}
\end{align}
where $k_{*}(t):=\argmax_{k \in \{1,...,\numK\}}  \| x(t)- a(t) \defGrad(t,\sigma^{(k)}(t))  \|^{2}$
By virtue of the above inequality (\ref{eq:bounded_state_transition}), the following property in Appendix A.3 of\cite{azuma2013broadcast} can be applied to (\ref{eq:bounded_state_transition}).
\begin{lem}[Boundedness of the state \cite{azuma2013broadcast}]\label{thm:convergence_of_azuma13}
	Let us define $t_{s}:=2s$ for $s \in \{0,1,...\}$. 
	Suppose that the conditions (c1), (c2), and (c3) hold.
	If the following relation
	\begin{equation}\label{eq:bounded_state_transition_BC}
	\begin{aligned}
	\|x(t_{s+1})\|^{2} {\leq} \| x(t_{s})- a(t_{s}) \defGrad(t_{s}, \sigma(t_{s}))  \|^{2}
	,
	\end{aligned}
	\end{equation}	
	is satisfied, then $\sup_{t_{s}} \|x(t_{s})\| < \infty$ surely holds, where the pair of the condition (c3) and the definition $t_{s}:=2s$ can be replaced by the pair of the condition (c3') and $t_{s}:=s$. 
\end{lem}	

Lemma \ref{thm:convergence_of_azuma13} proves the condition (c5) in Lemma \ref{thm:convergence_proof_of_BC} by regarding $t$, $t+1$, and $\sigma^{(k_{*}(t))}(t)$ in (\ref{eq:bounded_state_transition}) as $t_{s}$, $t_{s+1}$, and $\sigma(t_{s})$ in (\ref{eq:bounded_state_transition_BC}).
The condition (c6) holds in a manner similar to the case of the BC law (see Appendix A.2 of \cite{azuma2013broadcast} or Proposition 1 of \cite{Spall92}).
Lemma \ref{thm:convergence_proof_of_BC} is thus applied to the PBC law.
This completes the proof.
\qed

\section{Proof of Theorem \ref{thm:PBC_cost_improvement}} \label{pf:PBC_cost_improvement}
This theorem is proven by mathematical induction.
In the case of $t=0$, (\ref{eq:x_PBC_is_BC}) and (\ref{eq:J_PBC_is_BC}) clearly hold.
With (\ref{eq:input_PBC}), $x(2(t+1))$ by employing the BC law is given as follows \cite{azuma2013broadcast}:
\begin{equation}\label{eq:x_BC_transition}
\begin{aligned}
x(2(t+1))|_{\mathrm{BC}}
=x(2t)|_{\mathrm{BC}}
-a(2t)\defGrad(2t,\sigma(2t))|_{\mathrm{BC}}
.
\end{aligned}
\end{equation}
Using (\ref{eq:transition_PBC}), $x(t+1)$ under the PBC law is obtained as
\begin{equation}\label{eq:x_PBC_transition}
\begin{aligned}
x(t+1)|_{\mathrm{PBC}}
=x(t)|_{\mathrm{PBC}}
-a(t)\defGrad(t,\sigma^{(1)}(t))|_{\mathrm{PBC}}
.
\end{aligned}
\end{equation}
Since $a(t)|_{\mathrm{PBC}}=a(2t)|_{\mathrm{BC}}$,  $c(t)|_{\mathrm{PBC}}=c(2t)|_{\mathrm{BC}}$, and $\sigma^{(1)}(t)|_{\mathrm{PBC}}=\sigma(2t)|_{\mathrm{BC}}$, if (\ref{eq:x_PBC_is_BC}) and (\ref{eq:J_PBC_is_BC}) hold for $t$, they are also satisfied for $t+1$.
Therefore, (\ref{eq:x_PBC_is_BC}) and (\ref{eq:J_PBC_is_BC}) hold for all $t\geq 0$.
This completes the proof.
\qed

\section{Proof of Theorem \ref{thm:PBC_dist_reduction}} \label{pf:PBC_dist_reduction}
Let us define the following function in order to simplify the notation:
\begin{equation}\label{eq:def_positive_diff_prob}
\begin{aligned}
\Delta J(t,\sigma):=J(x(t)+c(t)\sigma) - J(x(t))
,
\end{aligned}
\end{equation}		
Since $\sigma_{i}=\sigma_{i}^{(-1)}$ and $\|\sigma_{i}\|=\sqrt{n}$ hold for any $i$, $\totaldist(2t)|_{\mathrm{BC}}$ is given by
\begin{equation}\label{eq:dist_BC}
\begin{aligned}
&
\totaldist(2t)|_{\mathrm{BC}}
\\&
=\sum_{\tau=0}^{t-1} \sum_{i=1}^{N}  \| c(2\tau)\sigma_{i}(2\tau) \| \Big|_{\mathrm{BC}}
\\&\;\;
+\Big\|c(2\tau)\sigma_{i}(2\tau)+a(2\tau)\frac{ \Delta J(2\tau,\sigma(2\tau))   }{c(2\tau)}\sigma_{i}^{(-1)}(2\tau) \Big\| \Big|_{\mathrm{BC}}
\\&
=\sum_{\tau=0}^{t-1} \sqrt{n} N \Big( c(2\tau) 
\\&\quad
+\Big|c(2\tau)+a(2\tau)\frac{ \Delta J(2\tau,\sigma(2\tau))   }{c(2\tau)}\Big|  \Big)
\Big|_{\mathrm{BC}}
.
\end{aligned}
\end{equation}

Here, for any positive scalar coefficients $\omega_{1}>0$ and $\omega_{2}>0$, the following expectation is derived:
\begin{equation}  \label{eq:prob_ineq}
\begin{aligned}
&
\mathrm{E}[ | \omega_{1} \Delta J(2\tau,\sigma) + \omega_{2}| ]
\\&=
\mathrm{Pr}\{\Delta J(2\tau,\sigma) \geq 0\}
\mathrm{E}[ | \omega_{1} \Delta J(2\tau,\sigma) + \omega_{2}| |_{\Delta J(2\tau,\sigma) \geq 0} ]
\\&\;
+
\mathrm{Pr}\{\Delta J(2\tau,\sigma) < 0\}
\mathrm{E}[ | \omega_{1} \Delta J(2\tau,\sigma) + \omega_{2}| |_{\Delta J(2\tau,\sigma) < 0} ]
\\&\geq
\mathrm{Pr}\{\Delta J(2\tau,\sigma) \geq 0\}
\mathrm{E}[  \omega_{1} \Delta J(2\tau,\sigma) + \omega_{2} |_{\Delta J(2\tau,\sigma) \geq 0} ]
\\&\;\;
-
\mathrm{Pr}\{\Delta J(2\tau,\sigma) < 0\}
\mathrm{E}[  \omega_{1} \Delta J(2\tau,\sigma) + \omega_{2} |_{\Delta J(2\tau,\sigma) < 0} ]
\\&=
\mathrm{E}[ |  \omega_{1} \Delta J(2\tau,\sigma) | ]
\\&\quad
+
(\mathrm{Pr}\{\Delta J(2\tau,\sigma) \geq 0\}-\mathrm{Pr}\{\Delta J(2\tau,\sigma) < 0\})
\omega_{2}
\\&=
\mathrm{E}[ |  \omega_{1} \Delta J(2\tau,\sigma) | ]
+
(2 \mathrm{Pr}\{\Delta J(2\tau,\sigma) \geq 0\} - 1)\omega_{2}
,
\end{aligned}
\end{equation}
where note that $| \omega_{1} \Delta J(2\tau,\sigma) + \omega_{2}| \geq -(\omega_{1} \Delta J(2\tau,\sigma) + \omega_{2})$ and $\mathrm{Pr}\{\Delta J(2\tau,\sigma) < 0\}=1-\mathrm{Pr}\{\Delta J(2\tau,\sigma) \geq 0\}$ hold.
In order to clarify the probability $\mathrm{Pr}\{\Delta J(2\tau,\sigma) \geq 0\}$ in (\ref{eq:prob_ineq}), we show the relation
\begin{equation} \label{eq:Prob_deltaJ}
\begin{aligned}
\mathrm{Pr}\{\Delta J(2\tau,\sigma(2\tau))|_{\mathrm{BC}} \geq 0\}
\geq \frac{\mathcal{P}(2\tau)}{2}
,
\end{aligned}
\end{equation}
where $\Delta J(2\tau,\sigma(2\tau))=J(x(2\tau)+c(2\tau)\sigma(2\tau)) - J(x(2\tau))$ in (\ref{eq:def_positive_diff_prob}).
If a random vector ${\sigma}$ is included in the set $\{-1,1\}^{nN}$, $-{\sigma}$ is also included in $\{-1,1\}^{nN}$ because each component of ${\sigma}$ is $-1$ or $1$.
If we assume that $\Delta J(2\tau,{\sigma})<0$ and $\Delta J(2\tau,-{\sigma})<0$ hold, then the relation 
\begin{equation} \label{eq:contradiction_of_quasi_conv}
\begin{aligned}
&
\max\{  J(x(2\tau)+c(2\tau){\sigma}) ,  J(x(2\tau)-c(2\tau){\sigma}) \}
\\&
<
J(x(2\tau))
\\&
=  J( \frac{x(2\tau)+c(2\tau){\sigma}}{2}+ \frac{x(2\tau)-c(2\tau){\sigma}}{2} )
\end{aligned}
\end{equation}
is satisfied.
This inequality (\ref{eq:contradiction_of_quasi_conv}) contradicts the necessary condition for the (quasi-)convexity of $J(x)$ on the minimum convex set including all possible $x(2\tau+1)|_{\mathrm{BC}}$.
Thus, for any $x(2\tau)|_{\mathrm{BC}}$ and any ${\sigma} \in \{-1,1\}^{nN}$, if $J(x)$ is convex or quasi-convex and $\Delta J(2\tau,{\sigma})|_{\mathrm{BC}}<0$ holds, then $\Delta J(2\tau,-{\sigma})|_{\mathrm{BC}}\geq 0$ holds.
Since the probability of ${\sigma}$ is constant on $\{-1,1\}^{nN}$, $\mathrm{Pr}\{\Delta J(2\tau,\sigma(2\tau))|_{\mathrm{BC}} \geq 0\} \geq \mathrm{Pr}\{\Delta J(2\tau,\sigma(2\tau))|_{\mathrm{BC}} < 0\}$, i.e., $\mathrm{Pr}\{\Delta J(2\tau,\sigma(2\tau))|_{\mathrm{BC}} \geq 0\} \geq 1/2$ holds if $J(x)$ is locally convex/quasi-convex, which yields (\ref{eq:Prob_deltaJ}).

Substituting the relation (\ref{eq:prob_ineq}) with $\omega_{1}=a(2\tau)/c(2\tau)$, $\omega_{2}=c(2\tau)$, and (\ref{eq:Prob_deltaJ}) into the expectation of $\totaldist(2t)|_{\mathrm{BC}}$ in  (\ref{eq:dist_BC}) yields
\begin{equation}\label{eq:E_dist_BC} 
\begin{aligned}
&
\mathrm{E}[ \totaldist(2t)|_{\mathrm{BC}} ]
\\&
\geq
\sum_{\tau=0}^{t-1} \sqrt{n} N
\Big( c(2\tau) 
+ (  \mathcal{P}(2\tau) - 1 ) c(2\tau)
\\&\qquad\qquad\qquad
+
a(2\tau)\frac{ \mathrm{E}[ |\Delta J(2\tau,\sigma(2\tau))| ] }{c(2\tau)}
\Big) \Big|_{\mathrm{BC}}
\\&
=
\sqrt{n} N
\\&\quad
\times
\sum_{\tau=0}^{t-1}
\Big( 
\mathcal{P}(2\tau) c(2\tau) 
+ a(2\tau)\frac{ \mathrm{E}[ |\Delta J(2\tau,\sigma(2\tau))| ]  }{c(2\tau)}
\Big)
\Big|_{\mathrm{BC}}
.
\end{aligned}
\end{equation}	
On the other hand,
\begin{align}
&
\totaldist(t)|_{\mathrm{PBC}}
\nonumber\\
&=\sum_{\tau=0}^{t-1} \sum_{i=1}^{N}  \|\sigma_{i}^{(1)}(\tau)\|
\Big|-a(\tau)\frac{ \Delta J(\tau,\sigma^{(1)}(\tau))   }{c(\tau)} \Big| \Big|_{\mathrm{PBC}}
\nonumber\\
&= \sqrt{n} N \sum_{\tau=0}^{t-1}
\Big|-a(\tau)\frac{ \Delta J(\tau,\sigma^{(1)}(\tau))   }{c(\tau)} \Big| \Big|_{\mathrm{PBC}}
,\label{eq:dist_PBC}
\\&
\mathrm{E}[ \totaldist(t)|_{\mathrm{PBC}} ]
=\sqrt{n} N \sum_{\tau=0}^{t-1}
a(\tau)\frac{ \mathrm{E}[ |\Delta J(\tau,\sigma^{(1)}(\tau))| ]  }{c(\tau)}
\Big|_{\mathrm{PBC}}
,\label{eq:E_dist_PBC}
\end{align}	
Recall the relation $x(t)|_{\mathrm{PBC}}=x(2t)|_{\mathrm{BC}}$ in (\ref{eq:x_PBC_is_BC}) and the assumptions that $a(t)|_{\mathrm{PBC}}=a(2t)|_{\mathrm{BC}}$,  $c(t)|_{\mathrm{PBC}}=c(2t)|_{\mathrm{BC}}$, and $\sigma^{(1)}(t)|_{\mathrm{PBC}}=\sigma(2t)|_{\mathrm{BC}}$ for all $t \geq 0$.  
From (\ref{eq:E_dist_BC}) and (\ref{eq:E_dist_PBC}), $\mathrm{E}[\totaldist(2t)|_{\mathrm{BC}}- \totaldist(t)|_{\mathrm{PBC}}]$ satisfies (\ref{eq:E_dist_PBCvsBC}).
The relation (\ref{eq:dist_PBCvsBC}) is proven because of (\ref{eq:dist_BC}) and (\ref{eq:dist_PBC}). 
Next, the following relation 
\begin{align}
&
\exists \tau \in \{0,..., t-1\}
,\;\mathrm{s.t.}\;
\Delta J(2\tau,\sigma(2\tau))|_{\mathrm{BC}} \geq 0
\nonumber\\&
\Rightarrow
\exists \tau \in \{0,..., t-1\}
,\;\mathrm{s.t.}\;
\nonumber\\&\quad
c(2\tau)|_{\mathrm{BC}} +\Big|c(2\tau)+a(2\tau)\frac{ \Delta J(2\tau,\sigma(2\tau))   }{c(2\tau)}\Big| \Big|_{\mathrm{BC}}
\nonumber\\&\quad
=2c(2\tau)|_{\mathrm{BC}} + \Big|a(2\tau)\frac{ \Delta J(2\tau,\sigma(2\tau))   }{c(2\tau)}\Big| \Big|_{\mathrm{BC}}
\nonumber\\&
\Rightarrow
\totaldist(2t)|_{\mathrm{BC}}
> \totaldist(t)|_{\mathrm{PBC}}
\label{eq:dist_BC_special2}
\end{align}
holds for the terms in (\ref{eq:dist_BC}). 
Therefore, $\Delta J(2\tau,\sigma(2\tau))|_{\mathrm{BC}}$ $<$ $0$, $\forall \tau \in \{0,..., t-1\}$ is the necessary condition for $\totaldist(2t)|_{\mathrm{BC}} = \totaldist(t)|_{\mathrm{PBC}}$.
The relation (\ref{eq:dist_PBCvsBC_prob}) is proven because $p(\Delta J(2\tau,\sigma(2\tau))|_{\mathrm{BC}} < 0) \leq 1 - ({\mathcal{P}(2\tau)}/{2})$ holds by transforming (\ref{eq:Prob_deltaJ}).
This completes the proof. 
\qed

\section{Proof of Theorem \ref{thm:MPBC_per_step_improvement}} \label{pf:MPBC_per_step_improvement}
We first prove the statement (i) for the case in which $1 \leq \numK_{b} < \numK_{a} < \infty$.
For brevity of notation, let us define $J_{\defGrad}( {\defGrad} ):=J(x(t+1))|_{x(t+1)=x(t)- a(t){\defGrad}}$ for a given $x(t)$, %
$\defGrad^{(k)}:=\defGrad(t,\sigma^{(k)}(t))$ in (\ref{eq:input_PBC}), and a combination ${\mathcal{C}}:={{}_{\numK_{a}}C_{\numK_{b}}}$.
For the given set $\{1,...,\numK_{a}\}$ of the index $k$ in $\defGrad^{(k)}$, we consider the ${\mathcal{C}}$ combinations to select $\numK_{b}$ indexes.
Let us define a vector $\mathcal{I}_{i} \in \{1,...,\numK_{a}\}^{\numK_{b}}$ as the array of the selected indexes in the $i$-th combination  
\begin{equation}
\begin{aligned}
&
\mathcal{I}_{i}
:=[\mathcal{I}_{i,1},
...,\mathcal{I}_{i,\numK_{b}}]^{\mathrm{T}}			
,\;
(i=1,2,...,{\mathcal{C}})
.
\end{aligned}
\end{equation}
Each index $k \in \{1,...,\numK_{a}\}$ is evenly selected $\numK_{b}{\mathcal{C}}/\numK_{a}$ times because the total number of selected indexes is $\numK_{b}{\mathcal{C}}$.
Thus, dividing $\defGrad^{(k)}$ by $\numK_{b}{\mathcal{C}}/\numK_{a}$ yields the following relation%
\footnote{%
	For example, if $\numK_{a}=4$ and $\numK_{b}=3$ are considered, ${\mathcal{C}}=4$, $\mathcal{I}_{i} \in \{[1,2,3]^{\mathrm{T}},$ $[1,2,4]^{\mathrm{T}},$ $[1,3,4]^{\mathrm{T}},$ $[2,3,4]^{\mathrm{T}}\}$, and 
	$
	\frac{1}{4}(\defGrad^{(1)}+\defGrad^{(2)}+\defGrad^{(3)}+\defGrad^{(4)})
	=\frac{1}{4}
	\{
	$ $
	\frac{\defGrad^{(1)}+\defGrad^{(2)}+\defGrad^{(3)}}{3}
	+$ $
	\frac{\defGrad^{(1)}+\defGrad^{(2)}+\defGrad^{(4)}}{3}
	+$ $
	\frac{\defGrad^{(1)}+\defGrad^{(3)}+\defGrad^{(4)}}{3}
	+$ $
	\frac{\defGrad^{(2)}+\defGrad^{(3)}+\defGrad^{(4)}}{3}
	\}$ hold.
}
\begin{equation}\label{eq:dividing_gradient}
\begin{aligned}
&
\frac{1}{\numK_{a}}\sum_{k=1}^{\numK_{a}}\defGrad^{(k)}
=
\frac{1}{\numK_{a}}
\sum_{i=1}^{\mathcal{C}} \sum_{j=1}^{\numK_{b}} 
\frac{\numK_{a} \defGrad^{(\mathcal{I}_{i,j})}}{\numK_{b} {\mathcal{C}}    }	
=
\sum_{i=1}^{\mathcal{C}} 
\frac{1}{\mathcal{C}}
\sum_{j=1}^{\numK_{b}} \frac{ \defGrad^{(\mathcal{I}_{i,j})} }{\numK_{b}}
.
\end{aligned}
\end{equation}	
Note that the function $J_{\defGrad}({\defGrad})$ is (strictly) convex in ${\defGrad}$ because $J(x(t+1))$ is (strictly) convex in $x(t+1)$ and ${\defGrad}=(x(t)-x(t+1))/a(t)$ is uniquely determined as the linear function of $x(t+1)$.
Since $\sum_{i=1}^{\mathcal{C}}({1}/{\mathcal{C}})=1$, applying Jensen's inequality to convex $J_{\defGrad}( {\defGrad})$ yields
\begin{equation}\label{eq:Jensen_convex}
\begin{aligned}
J_{\defGrad}\Big(
\sum_{i=1}^{\mathcal{C}} 
\frac{1}{\mathcal{C}}
\sum_{j=1}^{\numK_{b}} \frac{ \defGrad^{(\mathcal{I}_{i,j})} }{\numK_{b}}	
\Big)
\leq
\sum_{i=1}^{\mathcal{C}} 
\frac{1}{\mathcal{C}}
J_{\defGrad}\Big(	
\sum_{j=1}^{\numK_{b}} \frac{ \defGrad^{(\mathcal{I}_{i,j})} }{\numK_{b}}
\Big)		
.
\end{aligned}
\end{equation}	
	Here, we consider the case in which $J_{\defGrad}({\defGrad})$ is strictly convex.
	The equality in (\ref{eq:Jensen_convex}) holds only if $\sum_{j=1}^{\numK_{b}} { \defGrad^{(\mathcal{I}_{i,j})} }/{\numK_{b}}$ is constant for all $i$, i.e., $\defGrad^{(k)}$ is constant for all $k$.
	Because of the strict convexity, $\defGrad^{(k)} = \defGrad(t,\sigma^{(k)}(t)) \neq 0$ in (\ref{eq:input_PBC}) and the relation 
	$\defGrad^{(k)}=\defGrad^{(k^{\prime})} \Leftrightarrow \sigma^{(k)}=\sigma^{(k^{\prime})}$ hold.	
	Thus, if $J_{\defGrad}({\defGrad})$ is strictly convex, taking the expectation in (\ref{eq:Jensen_convex}) with respect to $\sigma^{(k)}$ (the set of $\sigma^{(k)}$ for each $k$ is finite) excludes the equality from (\ref{eq:Jensen_convex}):
	\begin{equation}\label{eq:Jensen_convex_expectation}
	\begin{aligned}
	&
	\mathrm{E}\Big[J_{\defGrad}\Big(
	\sum_{i=1}^{\mathcal{C}} 
	\frac{1}{\mathcal{C}}
	\sum_{j=1}^{\numK_{b}} \frac{ \defGrad^{(\mathcal{I}_{i,j})} }{\numK_{b}}	
	\Big)
	\Big|_{x(t)}
	\Big]	
	\\&
	{<}
	\mathrm{E}\Big[\sum_{i=1}^{\mathcal{C}} 
	\frac{1}{\mathcal{C}}
	J_{\defGrad}\Big(	
	\sum_{j=1}^{\numK_{b}} \frac{ \defGrad^{(\mathcal{I}_{i,j})} }{\numK_{b}}
	\Big)
	\Big|_{x(t)}	
	\Big]	
	.
	\end{aligned}
	\end{equation}	
	Using (\ref{eq:dividing_gradient}) and (\ref{eq:Jensen_convex_expectation}) yields the following inequality:
\begin{equation} \label{eq:E_Jensen_convex}
\begin{aligned}
&
\mathrm{E}[J(x(t+1)  |_{x(t)} ] |_{\numK=\numK_{a}}
=
\mathrm{E}\Big[
J_{\defGrad} \Big(\frac{1}{\numK_{a}}\sum_{k=1}^{\numK_{a}}\defGrad^{(k)} \Big)
\Big|_{x(t)}
\Big] 
\\&
=
\mathrm{E}\Big[
J_{\defGrad}\Big(	
\sum_{i=1}^{\mathcal{C}} 
\frac{1}{\mathcal{C}}
\sum_{j=1}^{\numK_{b}} \frac{ \defGrad^{(\mathcal{I}_{i,j})} }{\numK_{b}}	
\Big)
\Big|_{x(t)}
\Big] 
\\&
{\leq}
\mathrm{E}\Big[
\sum_{i=1}^{\mathcal{C}} 
\frac{1}{\mathcal{C}}
J_{\defGrad}\Big(	
 \sum_{j=1}^{\numK_{b}} \frac{ \defGrad^{(\mathcal{I}_{i,j})}	 }{\numK_{b}}
\Big)
\Big|_{x(t)}
\Big] 
\\&
=
\frac{1}{\mathcal{C}}
\sum_{i=1}^{\mathcal{C}} 
\mathrm{E}\Big[
J_{\defGrad}\Big(	
\sum_{k=1}^{\numK_{b}} \frac{\defGrad^{(k)}}{\numK_{b}}
\Big)
\Big|_{x(t)}
\Big] 
\\&
=
\mathrm{E}\Big[
J_{\defGrad}\Big(	
\frac{1}{\numK_{b}} \sum_{k=1}^{\numK_{b}} \defGrad^{(k)}
\Big)
\Big|_{x(t)}
\Big] 
=
\mathrm{E}[J(x(t+1) |_{x(t)} ] |_{\numK=\numK_{b}}
,
\end{aligned}
\end{equation}	
where the equality of $\leq$ in (\ref{eq:E_Jensen_convex}) is possible only if $J(x)$ is not strictly convex but convex.
If $J_{\defGrad}( {\defGrad} )$ is concave in ${\defGrad}$,  the inequality sign in (\ref{eq:E_Jensen_convex}) is reversed.
The statement (i) therefore holds.

 We next prove the statement (ii). 
In the above proof of the statement (i), let us redefine %
$J_{\defGrad}( {\defGrad} )$ %
$:= \sum_{i=1}^{N} \|x_{i}(t+1)-x_{i}(t)\|^{\kappa} |_{x_{i}(t+1)=x_{i}(t)- a(t){\defGrad}_{i}}$ %
$=a(t)^{\kappa} \sum_{i=1}^{N}  \|{\defGrad}_{i}\|^{\kappa}$ %
for a given $x(t)$, where ${\defGrad}=:[{\defGrad}_{i}^{\mathrm{T}},..., {\defGrad}_{N}^{\mathrm{T}}]^{\mathrm{T}}$.
Since the norm $\|{\defGrad}_{i}\|$ is convex in ${\defGrad}_{i}$ \cite{Boyd04}, $J_{\defGrad}( {\defGrad} )$ is convex in ${\defGrad}$ for $\kappa = 1$ and is strictly convex in ${\defGrad}$ for $\kappa > 1$.
The inequality in (\ref{eq:E_Jensen_convex}) holds in a similar manner, which leads to the statement (ii). 
This completes the proof.	
\qed

\section{Proof of Theorem \ref{thm:MPBC_multi_step_improvement}} \label{pf:MPBC_multi_step_improvement}
We consider the case in which $1 \leq \numK_{b} < \numK_{a} < \infty$.
We prove this theorem by focusing on the dynamic programming of stochastic dynamical systems.
Let us define 
\begin{align}
\valFunc_{s}(x(s))&:= \mathrm{E}[ J(x(t)) |_{x(s)}]
,
\end{align}
where $\valFunc_{t}(x(t)) = J(x(t))$ holds. We obtain
\begin{equation}
\begin{aligned}
\valFunc_{t-2}(x(t-2)) 
&
= \mathrm{E}[ \mathrm{E}[  \valFunc_{t}(x(t)) |_{x(t-1)} ] |_{x(t-2)} ]
\\&
= \mathrm{E}[              \valFunc_{t-1}(x(t-1))          |_{x(t-2)} ]
,
\end{aligned}
\end{equation}
Iterating the above transformation yields the recurrence formula of $\valFunc_{s}(x(s))$
\begin{equation}
\begin{aligned}
\valFunc_{s}(x(s))= \mathrm{E}[ \valFunc_{s+1}(x(s+1)) |_{x(s)}]
.
\end{aligned}
\end{equation}
Note that for all $s \in \{0,...,t-1\}$, $\valFunc_{s+1}(x(s+1))=\mathrm{E}[ J(x(t)) |_{x(s+1)}]$ is convex in $x(s+1)$ because, for all $\sigma^{(k)}(m)$ $(m \in \{0,...,t-1\},\; k \in \{1,...,\numK\})$, $J_{s+1,t}(\tilde{x}):=J(x(t))|_{x(s+1)=\tilde{x}}$ is convex in $\tilde{x}$ and the expectation of convex $J_{s+1,t}(\tilde{x})$ is also convex \cite{Boyd04}.
Let us define $x(s+1,\numK_{a})|_{x(s)}:=x(s+1)|_{x(s), \numK=\numK_{a}}$ and $\valFunc_{s}(x(s),\numK_{a}):=\mathrm{E}[ J(x(t)) |_{x(s)}]|_{\numK=\numK_{a}}$.
Applying Theorem \ref{thm:MPBC_per_step_improvement} to the expectation with respect to $\sigma^{(k)}(s)$ yields
\begin{equation} \label{eq:ineq_V_s1}
\begin{aligned}
&	
\mathrm{E}[\valFunc_{s+1}(x(s+1,\numK_{a}),\numK) |_{x(s)} ] 
\\&
\leq
\mathrm{E}[\valFunc_{s+1}(x(s+1,\numK_{b}),\numK) |_{x(s)} ] 
,
\end{aligned}
\end{equation}	
for all $s \in \{0,...,t-1\}$.
Meanwhile, by applying Theorem \ref{thm:MPBC_per_step_improvement} to the expectation with respect to $\sigma^{(k)}(t-1)$, the following relation holds for all $x(t-1,\numK)$ 
\begin{equation}  \label{eq:ineq_V_s2_pre}
\begin{aligned}
&
\valFunc_{t-1}(x(t-1,\numK),\numK_{a})
=
\mathrm{E}[J(x(t,\numK_{a}))  |_{x(t-1,\numK)}] 
\\&
\leq
\mathrm{E}[J(x(t,\numK_{b}))  |_{x(t-1,\numK)}] 
=
\valFunc_{t-1}(x(t-1,\numK),\numK_{b})	
.
\end{aligned}
\end{equation}
Thus,
\begin{equation}  \label{eq:ineq_V_s2}
\begin{aligned}
&
\mathrm{E}[\valFunc_{t-1}(x(t-1,\numK),\numK_{a})  |_{x(t-2)}] 
\\&
\leq
\mathrm{E}[\valFunc_{t-1}(x(t-1,\numK),\numK_{b})  |_{x(t-2)}] 
,
\end{aligned}
\end{equation}
holds.
Using (\ref{eq:ineq_V_s1}) and (\ref{eq:ineq_V_s2}) provides
\begin{equation} \label{eq:ineq_V_s3}
\begin{aligned}
\valFunc_{t-2}(x(t-2),\numK_{a})
&
=
\mathrm{E}[\valFunc_{t-1}(x(t-1,\numK_{a}),\numK_{a})  |_{x(t-2)}] 
\\&
\leq
\mathrm{E}[\valFunc_{t-1}(x(t-1,\numK_{b}),\numK_{a})  |_{x(t-2)}] 
\\&
\leq
\mathrm{E}[\valFunc_{t-1}(x(t-1,\numK_{b}),\numK_{b})  |_{x(t-2)}] 
\\&
=
\valFunc_{t-2}(x(t-2),\numK_{b})
.
\end{aligned}
\end{equation}
Using (\ref{eq:ineq_V_s1}) and (\ref{eq:ineq_V_s3}) yields	
\begin{equation} \label{eq:ineq_V_s4}
\begin{aligned}
\valFunc_{t-3}(x(t-3),\numK_{a})
&
=
\mathrm{E}[\valFunc_{t-2}(x(t-2,\numK_{a}),\numK_{a})  |_{x(t-3)}] 
\\&
\leq
\mathrm{E}[\valFunc_{t-2}(x(t-2,\numK_{b}),\numK_{a})  |_{x(t-3)}] 
\\&
\leq
\mathrm{E}[\valFunc_{t-2}(x(t-2,\numK_{b}),\numK_{b})  |_{x(t-3)}] 
\\&
=
\valFunc_{t-3}(x(t-3),\numK_{b})
.
\end{aligned}
\end{equation}	
Iterating the above process from $t-1$ to $0$ yields
\begin{equation} \label{eq:ineq_V_s5}
\begin{aligned}
&
\valFunc_{0}(x(0),\numK_{a}) 
= \mathrm{E}[ J(x(t))  |_{x(0)} ]|_{\numK=\numK_{a}}
\\&
\leq
\valFunc_{0}(x(0),\numK_{b}) 
= \mathrm{E}[ J(x(t))  |_{x(0)} ]|_{\numK=\numK_{b}}
.
\end{aligned}
\end{equation}
Since Theorem \ref{thm:MPBC_per_step_improvement} is used, the equalities in (\ref{eq:ineq_V_s1}), (\ref{eq:ineq_V_s2_pre}), (\ref{eq:ineq_V_s2}), (\ref{eq:ineq_V_s3}), (\ref{eq:ineq_V_s4}), and (\ref{eq:ineq_V_s5}) are possible only if $J_{s,t}(\tilde{x}):=J(x(t))|_{x(s)=\tilde{x}}$ is not strictly convex but convex in $\tilde{x}$.
Therefore, the relation (\ref{eq:ineq_V_s5}) is equivalent to (\ref{eq:multisteps_cost_imporvvement_by_MPBC}).
The relation (\ref{eq:multisteps_distance_imporvvement_by_MPBC}) can be proven as described above by applying $\totaldist(t)$ and the conditional expectations $\mathrm{E}[\cdot|_{x(0),...,x(s)} ]$ rather than $J(x(t))$ and $\mathrm{E}[\cdot|_{x(s)} ]$, respectively.
In this case, $\valFunc_{s+1}(x(s+1)):=\mathrm{E}[ \totaldist(t) |_{x(0),...,x(s+1)}]$ is indeed convex in $x(s+1)$ because $\totaldist_{s+1,t}(\tilde{x}):=(\totaldist(t)-\totaldist(s+1))|_{x(s+1)=\tilde{x}}$ and $\totaldist(s+1)|_{x(s+1)=\tilde{x}}$ for each fixed $x(0), ...,x(s)$ are convex in $\tilde{x}$.
This completes the proof.
\qed

\section{Proof of Proposition \ref{thm:quad_J}} \label{pf:quad_J}

The input $u_{i}(s)$ including $\nu^{(k)}(s)$ is linear in $x(s)$ because $J(x)$ is quadratic in $x$. 
Therefore, for all $s \in \{0,...,t\}$, $s^{\prime}  \in \{0,...,t\}$, and $\sigma^{(k)}(m)$ $(m \in \{0,...,t-1\},\; k \in \{1,...,\numK\})$, $x(s)$ is linear in $x(s^{\prime})$.
Since $J(x(t))$ is convex in $x(t)$, $J_{s,t}(\tilde{x})$ is convex in $\tilde{x}$.
Since $D(t)$ is convex in $x(s^{\prime})$ for all $s^{\prime} \in \{0,...,t\}$, $\totaldist_{s,t}(\tilde{x})$ is convex in $\tilde{x}$.
This completes the proof.	
\qed

\end{document}